\documentclass[a4paper]{article}
\usepackage[top=25truemm,bottom=25truemm,left=20truemm,right=20truemm,truedimen]{geometry}
\usepackage{amsmath,amsthm,amsfonts,amssymb,stmaryrd,indentfirst}
\usepackage[mathscr]{euscript}
\usepackage[shortlabels]{enumitem}
\usepackage{titlesec}
\usepackage{tikz-cd}
\usepackage{url}
\usepackage{hyperref}
\hypersetup{
    colorlinks = true,
    linkcolor = [rgb]{0.3,0.55,0.25},
    citecolor = [rgb]{0.164,0.39,0.77},
    urlcolor = [rgb]{0.7,0.2,0.4}
}

\DeclareMathOperator{\id}{id}
\DeclareMathOperator{\Ker}{Ker}
\DeclareMathOperator{\Imag}{Im}
\DeclareMathOperator{\Der}{Der}
\DeclareMathOperator{\DDer}{DDer}
\DeclareMathOperator{\Hom}{Hom}
\DeclareMathOperator{\IHom}{IHom}
\DeclareMathOperator{\End}{End}

\DeclareMathOperator{\Tr}{Tr}
\DeclareMathOperator{\Div}{Div}
\DeclareMathOperator{\sdiv}{div} 

\DeclareMathOperator{\ad}{ad}

\theoremstyle{plain}
\newtheorem{theorem}{Theorem}
\newtheorem*{theorem*}{Theorem}
\newtheorem{proposition}[theorem]{Proposition}
\newtheorem{deflem}[theorem]{Definition-Lemma}
\newtheorem{lemma}[theorem]{Lemma}
\newtheorem{corollary}[theorem]{Corollary}

\numberwithin{theorem}{section}

\theoremstyle{definition}
\newtheorem{definition}[theorem]{Definition}
\newtheorem{remark}[theorem]{Remark}
\newtheorem{example}[theorem]{Example}

\everymath{\displaystyle}
\allowdisplaybreaks[0]

\titleformat{\section}{\large\scshape}{\thesection.}{3pt}{}

\usepackage{titling}  
\pretitle{
	\begin{center}\Large\bfseries
}
\posttitle{
	\end{center}\ \\[-20pt]
}
\preauthor{
	\begin{center}\large
}
\postauthor{
	\end{center}\ \\[-60pt]
}

\title{Non-commutative Divergence and the Turaev Cobracket}
\author{Toyo TANIGUCHI \thanks{Graduate School of Mathematical Sciences, The University of Tokyo. 3-8-1, Komaba, Meguro-ku, Tokyo, 153-8914, Japan. E-mail: \texttt{toyo(at)ms.u-tokyo.ac.jp}}}
\date{}

\begin{document}
\maketitle

\begin{abstract}
\noindent The divergence map, an important ingredient in the algebraic description of the Turaev cobracket on a connected oriented compact surface with boundary, is reformulated in the context of non-commutative geometry using a flat connection on the space of $1$-forms on a formally smooth associative algebra. We then extend this construction to homologically smooth associative algebras, which allows us to give a similar algebraic description of the Turaev cobracket on a closed surface. We also look into a relation between the Satoh trace and the divergence map on a free Lie algebra via geometry over the Lie operad.\\
\end{abstract}


\section{Introduction}

\indent The Turaev cobracket is a \textit{loop operation} introduced in \cite{Turaev:1991} by Turaev, which endows a structure of an involutive Lie bialgebra on the space of non-trivial free loops on an oriented surface together with the Goldman bracket (the involutivity is due to Chas \cite{Chas:aa}). These topologically defined operations have several interesting algebraic descriptions by Massuyeau \cite{Massuyeau:2018aa} and by Kawazumi and Kuno \cite{Kawazumi:2016}, a combinatorial description by Chas \cite{Chas:aa}, and also relate to the necklace Lie bialgebra in quiver theory.

Since connected oriented and closed surfaces are distinguished by their fundamental group, one would expect to recover a significant amount of topological information from the group, including those loop operations. In that case, an algebraic description of the Goldman bracket was given by Vaintrob \cite{Vaintrob:aa}, and it is the one induced from the commutator bracket of derivations combined with the Poincar\'e--Van den Bergh duality. In the case of connected oriented compact surfaces with non-empty boundaries, another algebraic description of the Turaev cobracket is given in the remarkable paper \cite{akkn} by Alekseev, Kawazumi, Kuno and Naef, where an interesting relation between two-dimensional topology and Kashiwara--Vergne groups is revealed. They introduced a non-commutative analogue of the divergence map, which is the main object of this paper. 

We now quickly recall the notion of a connection, or a covariant derivative, in terms of non-commutative geometry in the sense of Kontsevich and Ginzburg (and many others); for the detail, see Sections \ref{sec:ncgeom}, \ref{sec:condiv} and \ref{sec:perf} in the body. Let $B$ be a unital associative algebra over a field $\mathbb{K}$ of characteristic zero, $\Omega^1B$ the space of $1$-forms on $B$, and $M$ a left $B$-module. A \textit{connection} is a $\mathbb{K}$-linear map
\[
	\nabla\colon M \to \Omega^1B\otimes_B M
\]
satisfying the Leibniz rule $\nabla(bm) = db\otimes m + b\nabla(m)$. The associated divergence of a derivation $f:B\to B$ is defined by the formula
\[
	\Div^\nabla(f) = \Tr(L_f - (i_f\otimes\id)\circ\nabla)\,,
\]
in which we use the Lie derivative $L_f$, the contraction map $i_f$, and the trace map whenever defined; in fact, the trace is well-defined if and only if $M$ is dualisable (i.e. finitely generated and projective) over $B$. This map gives a Lie algebra $1$-cocycle if the connection is \textit{flat}: $\nabla^2=0$.

We then apply this to the case $B = A^\mathrm{e}$ and $M = \Omega^1A$ for another unital associative algebra $A$. In this setting, $\Omega^1A$ being dualisable is equivalent to the condition that $A$ is cohomologically $1$-dimensional, called \textit{formal smoothness}. We further extend this construction in the case of homologically smooth algebras by considering a projective resolution of $M$ and a family of connections on it, which we call a \textit{homological connection}, similar to that having finite cohomological dimension is called homologically smooth.

The divergence map in the context of geometry over an operad was already studied in \cite{Powell:2021aa} by Powell, where they defined the standard divergence map on free algebras and showed the $1$-cocycle property. This standard divergence can be recovered in our construction by taking a canonical flat connection on a free algebra; for the concrete formulation, see Section \ref{sec:flatconns}. Our motivation to introduce another formulation of divergences is to deal with other associative algebras, namely, the group algebra on free groups, which are isomorphic to the fundamental group of a surface with boundaries.

The Turaev cobracket and the divergence map are deeply related to the family of ``trace maps'' in the Johnson theory. The original trace map goes back to the paper \cite{Morita_1993} by Morita, where their motivation was to study the cokernel of the Johnson homomorphism in the context of the mapping class groups of a surface. Later, the Morita trace is refined in their study of the automorphism group of a free group in the paper \cite{SATOH2012709} by Satoh and applied to the mapping class group of a surface in \cite{Enomoto:2010aa} by Enomoto and Satoh.

Our construction of the divergence map leads us to an aspect of the Satoh trace on the automorphism group of a free group. It was combinatorially defined in their original paper \cite{SATOH2012709}, but it can be realised as the divergence map associated with a non-commutative version of a flat connection, briefly discussed in Section \ref{sec:Lie}. For the relation between the Turaev cobracket and the mapping class groups of surfaces, see Theorem 1.14 and Section 9 of \cite{akkn} for a detailed explanation. A pre-existing interpretation of the Satoh trace is studied in the paper \cite{Massuyeau:2020aa} by Massuyeau and Sakasai, where they introduce several variants of the trace map and relate it to the Magnus representation of the group of certain automorphisms on a free Lie algebra.

Finally, the goal of this paper is to give an algebraic description of the Turaev cobracket in the case of connected oriented and closed surfaces by utilising the tools developed above. Let $\mathbb{K}$ be a unital ring, $\Sigma$ such a surface, and $\pi=\pi_1(\Sigma)$ the fundamental group of a surface. The symbol $|\cdot|$ denotes the cyclic quotient so that the space $|\mathbb{K}\pi|$ is the free $\mathbb{K}$-module spanned by the homotopy classes of free loops on $\Sigma$. The main result is the following. Recall that the divergence takes the form $\Der_\mathbb{K}(\mathbb{K}\pi)\to |\mathbb{K}\pi^\mathrm{e}|\cong |\mathbb{K}\pi|^{\otimes 2}$, and the first Hochschild cohomology $\mathrm{HH}^1(\mathbb{K}\pi)$ of the group ring $\mathbb{K}\pi$ is regarded as the space of outer derivations.

\begin{theorem*}
We have a $\mathbb{K}$-linear map
\[
	\Div\colon \mathrm{HH}^1(\mathbb{K}\pi)\to |\mathbb{K}\pi/\mathbb{K}1|^{\otimes 2}
\]
constructed from a flat homological connection (see Definition \ref{def:homconn}), such that the composite
\[\begin{tikzcd}[ampersand replacement = \&]
	\text{\textbar}\mathbb{K}\pi\text{\textbar} \arrow[r, "v"] \& \mathrm{HH}^1(\mathbb{K}\pi) \arrow[r,"\Div"] \& \text{\textbar}\mathbb{K}\pi/\mathbb{K}1\text{\textbar}^{\otimes 2} .
\end{tikzcd}\]
with the Vaintrob's map $v$ recalled in Section \ref{sec:closed} is the Turaev cobracket.
\end{theorem*}
\noindent For more precise statement, see Theorem \ref{thm:cob} and Corollary \ref{cor:cob} in the body.\\

\noindent\textbf{Organisation of the paper.} In Section \ref{sec:cob}, we recall the Turaev cobracket and the divergence map defined in \cite{akkn}. Sections \ref{sec:ncgeom}-\ref{sec:flatconns} introduce the language of non-commutative geometry and reformulate the divergence map. In Section \ref{sec:Lie}, which is logically independent of later sections, we look into divergence maps on Hopf algebras and then in geometry over the Lie operad and its relation with the Satoh trace. Finally, the case of closed surfaces is dealt with in Sections \ref{sec:perf} and \ref{sec:closed}.\\

\noindent\textbf{Acknowledgements.} The author would like to thank Nariya Kawazumi for insightful advice, Florian Naef for useful comments on the draft of this paper, Yusuke Kuno for generously sharing their notes on this topic, Anton Alekseev for pointing out some important references, and the referee for improving many details and suggesting the re-organisation of Section \ref{sec:Lie}.\\

\noindent\textbf{Conventions.} $\mathbb{K}$ is a field of characteristic zero throughout this paper. All $\mathbb{K}$-algebras contain $\mathbb{K}1$ in their centre. Unadorned tensor products are always over $\mathbb{K}$.\\

\section{The Turaev Cobracket and its Algebraic Description}\label{sec:cob}
Let $\Sigma$ be a connected oriented surface possibly with boundary, and $\pi = \pi_1(\Sigma)$ its fundamental group. Put $|\mathbb{K}\pi| = \mathbb{K}\pi/[\mathbb{K}\pi,\mathbb{K}\pi]$, the cyclic quotient of the algebra $\mathbb{K}\pi$, which is a free $\mathbb{K}$-module spanned by the homotopy classes of free loops on $\Sigma$. The Turaev cobracket is a map $\delta\colon |\mathbb{K}\pi| \to |\mathbb{K}\pi/\mathbb{K}1|^{\otimes 2}$ defined by, for a generically immersed free loop $\alpha\colon [0,1]/\{0,1\}\to \Sigma$,
\begin{align*}
	\delta(\alpha) = \sum_{\substack{t_1\neq t_2\in[0,1]\\\alpha(t_1)=\alpha(t_2)}}  \mathrm{sign}(\alpha;t_1,t_2)\,\alpha|_{[t_1,t_2]}\otimes\alpha|_{[t_2,t_1]}\,,
\end{align*}
where $\mathrm{sign}(\alpha;t_1,t_2)$ is the local intersection number with respect to the orientation of $\Sigma$. This map is well-defined up to birth-deaths of monogons, hence takes its value in $|\mathbb{K}\pi/\mathbb{K}1|^{\otimes 2}$. If $\Sigma$ admits a framing $\mathsf{fr}$ (i.e., a smooth non-vanishing vector field), we can upgrade it to the map $$\delta^\mathsf{fr}\colon|\mathbb{K}\pi| \to |\mathbb{K}\pi|^{\otimes 2}$$ by taking a rotation-free representative of $\alpha$. 

Now let $n\geq 0$ and assume that $\Sigma$ is a connected oriented compact surface $\Sigma_{g,n+1}$ of genus $g$ with $(n+1)$ boundary components numbered through $0$ to $n$, and take a base point on the $0$-th boundary component. Then, an algebraic description of $\delta^\mathsf{fr}$ involving the non-commutative divergence is given in \cite{akkn} as follows. Firstly, $\sigma$ is defined as the based version of the Goldman bracket
\begin{align}\label{eq:sigma}
	\sigma\colon |\mathbb{K}\pi| \to \Der_\mathbb{K}(\mathbb{K}\pi)
\end{align}
from the space of free loops to the space of $\mathbb{K}$-linear derivations on $\mathbb{K}\pi$. It is given by 
\[
	\sigma(\alpha)(x) = \sum_{p\in\alpha\cap x} \mathrm{sign}(\alpha,x;p) \,\alpha\ast_px
\]
for generic representatives of a free loop $\alpha$ and $x\in\pi$. Here $\alpha\ast_px$ is a based loop obtained by traversing $x$ until it reaches $p$, then going along $\alpha$, and finally following the rest of $x$. (In \cite{akkn}, $\sigma$ is described as the map induced from the double bracket $\kappa$, but we omit the detail here.)

Next, take the free-generating system $\mathcal{C} = (a_i,b_i,\zeta_j)_{1\leq i\leq g,1\leq j\leq n}$ of $\pi$ so that $\zeta_j$ is homotopic to the loop represented by the $j$-th boundary component with the induced orientation, and
\[
	\zeta_0 = (a_1,b_1)\cdots(a_g,b_g)\zeta_1\cdots\zeta_n
\]
holds (see Figure 10 of \cite{akkn}). Here $(x,y) = xyx^{-1}y^{-1}$ is the group commutator. This induces an isomorphism $\pi\cong F_{2g+n}$, which we fix throughout this paper. For all $\alpha$, the derivation $\sigma(\alpha)$ vanishes on the subalgebra $\mathbb{K}\langle\zeta_0\rangle$ generated by $\zeta_0$ in $\mathbb{K}\pi$, which we write as 
\begin{align*}
	\sigma\colon |\mathbb{K}\pi| \to \Der_{\mathbb{K}\langle\zeta_0\rangle}(\mathbb{K}\pi)\,.
\end{align*}
The divergence map associated with $\mathcal{C}$ is defined by, for $f\in \Der_\mathbb{K}(\mathbb{K}\pi)$,
\begin{align}\label{eq:divc}
	\Div^\mathcal{C}(f) = \sum_{c\in \mathcal{C}}|\partial_c(f(c)) - 1\otimes c^{-1}f(c)|\in|\mathbb{K}\pi|\otimes|\mathbb{K}\pi|\,,
\end{align}
where $\partial_c$ is the double derivation (see the end of the next section) defined by the formula $\partial_c(c') = \delta_{c,c'}1\otimes 1$.

\begin{theorem}\label{thm:akkn}
\rm {(Theorem 5.16, \cite{akkn})} \it Let $\mathsf{fr}$ be the framing such that all generators $x_i$ are rotation-free. Then, the composite $\Div^\mathcal{C}\circ\,\sigma$ is equal to the framed Turaev cobracket $\delta^\mathsf{fr}$.\\
\end{theorem}

\section{Preliminaries on Non-commutative Geometry}\label{sec:ncgeom}
In this section, we recall some definitions in non-commutative geometry. Let $A$ be a unital associative $\mathbb{K}$-algebra with the multiplication map $\mu\colon A\otimes A\to A$, and $A^\mathrm{e} = A\otimes A^\mathrm{op}$ its enveloping algebra. We use the following notation: $\bar x$ is the copy of $x\in A$ in $A^\mathrm{op}$, and elements of $A^\mathrm{e}$ are written without the tensor symbol. Therefore, we have the relations $x\bar y = \bar y x$ and $\bar x\bar y = \overline{yx}$. We identify $A$-bimodules with \textit{left} $A^\mathrm{e}$-modules. 

\begin{definition} \ 
	\begin{itemize}
		\item The left $A^\mathrm{e}$-module structure on $A$ is given by
		\[
			x\bar y\cdot a = xay \quad\textrm{ for } a\in A\textrm{ and } x\bar y\in A^\mathrm{e}\,.
		\]
		\item The left $A^\mathrm{e}$-module structure on $A\otimes A$, the \textit{outer} structure, is given by
		\[
			x\bar y\cdot(a\otimes b) = xa\otimes by \quad\textrm{ for } a\otimes b\in A\otimes A\textrm{ and } x\bar y\in A^\mathrm{e},
		\]
		while the right $A^\mathrm{e}$-module structure, the \textit{inner} structure, is given by $(a\otimes b)\cdot x\bar y = ax\otimes yb$. The natural identification $A^\mathrm{e}\cong A\otimes A$ is an isomorphism of $A^\mathrm{e}$-bimodules.
		\item $\Omega^1A = \Ker(\mu\colon A\otimes A\to A)$ is the space of \textit{non-commutative $1$-forms}. This is a left $A^\mathrm{e}$-submodule of $A\otimes A$, since $\mu$ is a left $A^\mathrm{e}$-module homomorphism. The universal derivation $d\colon A\to \Omega^1A$ is defined by $da = 1\otimes a - a\otimes 1$.
		\item $\Omega^\bullet A = \bigoplus_{m\geq 0}(\Omega^1A)^{\otimes_Am}$ is the tensor algebra over $A$ generated by $1$-forms, which is a graded $A$-algebra.
	\end{itemize}
\end{definition}

Let $\bar A = A/\mathbb{K}1$. Then $\Omega^1A$ is canonically isomorphic to $A\otimes \bar A$ as a left $A$-module by the map given by 
\begin{align*}
	A\otimes \bar A\to \Omega^1A\colon a_0\otimes [a_1] &\mapsto a_0da_1 := a_0\otimes a_1 - a_0a_1\otimes 1\,.
\end{align*}
Its inverse is given by the projection $A\otimes A\twoheadrightarrow A\otimes\bar A$. Similarly, $\Omega^1A$ is also isomorphic to $\bar A\otimes A$ as a right $A$-module. With this notation, every element of $\Omega^\bullet A$ can be written as a linear combination of elements of the form $a_0da_1\dotsc da_n$, abbreviating the tensor symbol.\\

Now let $B$ be another associative $\mathbb{K}$-algebra, $M$ a left $B$-module, and $M^* = \Hom_B(M,B)$ the dual space of $M$. Then, $M^*$ is naturally a right $B$-module by
\[
	(\theta\cdot b)(m) = \theta(m)b \quad \textrm{ for } \theta\in M^*\!,\, b\in B \textrm{ and } m\in M.
\]

\begin{definition}
A left $B$-module $M$ is said to be \textit{dualisable} if it is finitely generated and projective over $B$.
\end{definition}

\begin{proposition}\label{prop:dualisable}
For a left $B$-module $M$, the following statements are equivalent:
\begin{enumerate}[(1)]
	\item $M$ is dualisable;
	\item There exists a non-negative integer $r$ and maps $e\colon B^{\oplus r}\to M$, $e^*\colon M\to B^{\oplus r}$ such that $e\circ e^* = \id_M$; and
	\item For any left $B$-module $N$, the map
	\[
		M^*\otimes_B N \to \Hom_B(M,N)\colon \theta\otimes n\mapsto (m\mapsto \theta(m)n)
	\]
	is an isomorphism.
\end{enumerate}
\end{proposition}
\noindent Proof. $(1)\Rightarrow(2)$: Since $M$ is finitely generated, there exists a non-negative integer $r$ and a surjective $B$-module homomorphism $e\colon M\to B^{\oplus r}$. By the projectivity of $M$, we can take a section $e^*\colon M\to B^{\oplus r}$ of $e$.\\
$(2)\Rightarrow(3)$ Denote by $e^i = e(0,\dotsc,1\dotsc,0)\in M$ the value on the $i$-th standard unit vector on $B^{\oplus r}$ and by $e^*_i\colon M\to B$ the projection of $e^*$ to the $i$-th component. The following map
\begin{align*}
	 \Hom_B(M,N)\to M^*\otimes_B N\colon \psi \to \sum_i e^*_i\otimes \psi(e_i)
\end{align*}
gives the inverse, which follows from the condition $e\circ e^* = \id_M$.\\
$(3)\Rightarrow(2)$: Write the corresponding element to $\id_M\in\Hom_B(M,M)$ in $M^*\otimes_BM$ by $\sum_{1\leq i\leq r} e^*_i\otimes e^i$. This defines the maps
\[	
	e\colon B^{\oplus r}\to M:(0,\dotsc,1\dotsc,0)\mapsto e^i\;\textrm{ and }\;e^*\colon M\to B^{\oplus r}:m \mapsto (e^*_i(m))_{1\leq i\leq r}
\]
which are what we want.\\
$(2)\Rightarrow(1)$: We have the following split exact sequence of $B$-module maps:
\[
	0\to \Ker e \to B^{\oplus r}\to M \to 0\,.
\]
Therefore $M$ is a direct summand of $B^{\oplus r}$.\qed\\

We will use the above proposition later to define the trace map.

\begin{definition}
An algebra $A$ is said to be \textit{formally smooth} if it is finitely generated as a $\mathbb{K}$-algebra with $A^\mathrm{e}$-projective $\Omega^1A$.
\end{definition}

This smoothness condition is convenient yet very restrictive. First, if $A$ is a finitely generated algebra, $\Omega^1A$ is finitely generated as an $A^\mathrm{e}$-bimodule by the Leibniz rule. Therefore, $A$ being formally smooth forces $\Omega^1A$ to be dualisable over $A^\mathrm{e}$. We will see two examples of formally smooth algebras.
\begin{example} 
Let $r\geq 0$. The free associative algebra $A = \mathbb{K}\langle z_1,\dotsc,z_r\rangle$ of rank $r$ is formally smooth. In fact, we have the standard $A^\mathrm{e}$-free resolution of $A$:
\begin{align*}
\begin{tikzcd}[cramped,ampersand replacement = \&]
	0 \arrow[r] \&  A\otimes\mathbb{K}\{z_1,\dotsc,z_r\}\otimes A \arrow[r,"\delta"] \& A\otimes A \arrow[r,"\mu"] \& A \arrow[r] \& 0\,,
\end{tikzcd}
\end{align*}
\[
	\delta(1\otimes z_i\otimes 1) = 1\otimes z_i - z_i\otimes 1 = dz_i\,,
\]
which shows that $\Omega^1A$ is an $A^\mathrm{e}$-free module with the basis $(dz_i)_{1\leq i\leq r}$.
\end{example}

The other example is the group algebra of a finitely generated free group. Before the proof, we define a useful functor.

\begin{definition}\label{def:functorphi}
Let $(A,\mu,\Delta,\eta,\varepsilon,S)$ be a Hopf algebra over $\mathbb{K}$. The functor $\Phi_A$ is defined by 
\begin{align}
	\Phi_A: A\textrm{-}\textbf{Mod}\to A^\mathrm{e}\textrm{-}\textbf{Mod}: M \mapsto \Phi_A(M),\;\; \psi\mapsto\psi\otimes\id_{A}\,.\;\label{eq:hopfadj}
\end{align}
Here we set $\Phi_A(M)=M\otimes A$ as a $\mathbb{K}$-module with the $A^\mathrm{e}$-action given by, for 	$a,x,y\in A$ and $m\in M$,
\[
	x\cdot(m\otimes a)\cdot y = x^{(1)}m\otimes x^{(2)}ay,
\]
using the coproduct $\Delta(x) = x^{(1)}\otimes x^{(2)}$.
\end{definition}

\begin{remark}\label{rem:adj}
We have canonical isomorphisms of $A$-bimodules: for $x,y\in A$ and $k\in \mathbb{K}$,
\begin{align}\label{eq:bimodisom}
	\begin{aligned}
		F\colon \Phi_A(A) \cong A^{\otimes 2}&:  x\otimes y \mapsto x^{(1)}\otimes S(x^{(2)})y\textrm{ and}\\
		\Phi_A(\mathbb{K}) \cong A\hspace{10pt}&:  k\otimes x  \mapsto kx\,.
	\end{aligned}
\end{align}
In addition, $\Phi_A$ admits the right adjoint $\Psi_A$ given, for an $A$-bimodule $N$, by $\Psi_A(N) = N$ with the conjugation action
\[
	x\cdot n = x^{(1)}nS(x^{(2)})\quad\textrm{for } x\in A\textrm{ and } n\in N
\]
and $\Psi_A(\psi) =\psi$ for a morphism $\psi$.\\
\end{remark}

\begin{proposition}\label{prop:freegroupsmooth}
Let $r\geq 0$. Then, the group algebra $A = \mathbb{K}F_r$ of the free group $F_r$ is formally smooth.
\end{proposition}
\noindent Proof. The following is the proof the author learned from Florian Naef. First, we fix a free-generating system $(x_i)_{1\leq i\leq r}$ of $F_r$. Then we have an exact sequence of left $A$-modules:
\begin{align*}
	\begin{tikzcd}[cramped,ampersand replacement = \&]
0 \arrow[r] \&  \bigoplus_{1\leq i\leq r}A\otimes\mathbb{K}\cdot(x_i-1) \arrow[r,"\iota"] \& A\arrow[r,"\varepsilon"] \& \mathbb{K} \arrow[r] \& 0\,,
	\end{tikzcd}
\end{align*}
where $\iota$ is the summation map and $\varepsilon$ is the augmentation map. The inverse of $\iota$ is given by the direct sum of Fox derivatives
\[
	\frac{\partial}{\partial x_i}:\Ker\varepsilon \to A\otimes\mathbb{K}\cdot(x_i-1)\,,
\]
which is uniquely specified by the formulae
\begin{align*}
	\frac{\partial}{\partial x_i}(ab) = \frac{\partial}{\partial x_i}(a) + a\frac{\partial}{\partial x_i}(b)\quad\textrm{and}\quad \frac{\partial}{\partial x_i}(x_j) = \delta_{ij}\,,
\end{align*}
and satisfies the equation
\[
	a - \varepsilon(a) = \sum_{1\leq i\leq r} \frac{\partial}{\partial x_i}(a)\cdot (x_i-1)\,.
\]
By applying the functor in Definition \ref{def:functorphi} to the resolution above, we obtain an acyclic complex
\begin{align*}
	\begin{tikzcd}[cramped,ampersand replacement = \&]
	0 \arrow[r] \& \Phi_A\left( \bigoplus_{1\leq i\leq r}A\cdot(x_i-1)\right) \arrow[r,"\iota\otimes \id"] \& \Phi_A(A)  \arrow[r,"\varepsilon\otimes\id"] \& \Phi_A(\mathbb{K}) \arrow[r] \& 0\,.
	\end{tikzcd}
\end{align*}
Transporting the differentials using the isomorphisms (\ref{eq:bimodisom}), we have a left $A^\mathrm{e}$-free resolution of $A$:
\begin{gather*}
	\begin{tikzcd}[cramped,ampersand replacement = \&]
	0 \arrow[r] \& \bigoplus_{1\leq i\leq r}A\otimes\mathbb{K}\cdot(x_i-1)\otimes A \arrow[r,"\delta'"] \& A^{\otimes 2}  \arrow[r,"\mu"] \& A \arrow[r] \& 0\,,
	\end{tikzcd}\\
	\delta'(1\otimes(x_i-1)\otimes 1) = x_i\otimes x_i^{-1} - 1\otimes 1 = -dx_i \,x_i^{-1}\,.
\end{gather*}
Therefore $\Omega^1A$ is an $A^\mathrm{e}$-free module with the basis $(dx_i \,x_i^{-1})_{1\leq i\leq r}$.\qed\\

We need some more definitions from non-commutative geometry.

\begin{definition} \ 
	\begin{itemize}
		\item For a $\mathbb{K}$-subalgebra $R$ of $A$ and an $A$-bimodule $M$, an $R$-linear derivation on $A$ into $M$ is a $\mathbb{K}$-linear map $f\colon A\to M$ such that
		\begin{align*}
			f(ab)& = f(a)\cdot b + a\cdot f(b)\quad\textrm{ for } a,b\in A,\textrm{ and}\\
			f(r) &= 0 \hspace{78pt} \textrm{ for } r\in R\,.
		\end{align*}
		The space of all such derivations is denoted by $\Der_R(A,M)$. Set $\Der_R(A) = \Der_R(A,A)$.
		\item For $f\in\Der_\mathbb{K}(A,M)$, the $A$-bimodule homomorphism $i_f\colon \Omega^1 A\to M$ is defined by
		\[
			i_f(da) = f(a)\quad\textrm{ for } a\in A\,,
		\]
		which gives a canonical isomorphism of $\mathbb{K}$-modules
		\[
			\Der_\mathbb{K}(A,M)\to\Hom_{A^\mathrm{e}}(\Omega^1A,M)\colon f\mapsto i_f
		\]
		with the inverse given by the composition with the universal derivation $d$.
		Conversely, the pair $(\Omega^1A, d)$ is characterised by the following universal property: for any $A$-bimodule $M$ and a $\mathbb{K}$-linear derivation $f\colon A\to M$, there is a unique $A^\mathrm{e}$-module map $i_f\colon \Omega^1A\to M$ such that $i_f\circ d = f$ holds.
		\item For $f\in\Der_\mathbb{K}(A)$, the corresponding map $i_f\colon \Omega^1A\to A$ is extended to the $A$-linear, degree $(-1)$ derivation $i_f\colon\Omega^\bullet A \to \Omega^\bullet A$ with Koszul signs, which is called the \textit{contraction map}.
		\item A \textit{double derivation} is a $\mathbb{K}$-linear map $\theta\colon A\to A\otimes A$ satisfying 
		\[
			\theta(ab) = \theta(a)\cdot b + a\cdot \theta(b) \quad\textrm{ for } a,b\in A.
		\]
		$\DDer_\mathbb{K}(A)$ denotes the space of $\mathbb{K}$-linear double derivations on $A$, which is isomorphic to $\Hom_{A^\mathrm{e}}(\Omega^1A,A^\mathrm{e})$. This is an $A$-bimodule by the inner structure of $A\otimes A$.
		\item The universal derivation $d\colon A \to \Omega^1A$ is extended to the $A$-linear, degree $1$ derivation $d\colon\Omega^\bullet A \to \Omega^\bullet A$ with Koszul signs, which is called the \textit{exterior derivative} on $A$.
		\item For $f\in\Der_\mathbb{K}(A)$, the \textit{Lie derivative} $L_f \colon \Omega^\bullet A\to \Omega^\bullet A$ is the $\mathbb{K}$-linear, degree $0$ derivation given by $L_f = [d,i_f] = d\circ i_f + i_f\circ d$, and it is determined by
		\[
			L_f(a) = f(a)\;\textrm{ and }\; L_f(da) = df(a)\,.
		\]
		Then, the map
		\[
			L\colon\Der_\mathbb{K}(A)\to\Der_\mathbb{K}(\Omega^\bullet A)^{(0)}\colon f\mapsto L_f
		\]
		is a Lie algebra homomorphism. The superscript $(0)$ denotes the degree $0$ part.
		\item The Hochschild homology and cohomology of $A$ with coeffitients in $A$ are defined by
		\[
			\mathrm{HH}_\bullet(A) = \mathrm{Tor}_\bullet^{A^\mathrm{e}}(A,A)\quad \textrm{and}\quad\mathrm{HH}^\bullet(A) = \mathrm{Ext}^\bullet_{A^\mathrm{e}}(A,A).\\
		\]
		The space $\mathrm{HH}_0(A)$ is called the \textit{trace space} of $A$ and is briefly denoted by $|A|$, regarded as the cyclic quotient. The space $\mathrm{HH}^1(A)$ is identified with the space of outer derivations.\\
	\end{itemize}
\end{definition}

\section{Connections and Divergences}\label{sec:condiv}
In this section, we recall the definition of a non-commutative connection, which we use to define the associated divergence map further. Let $B$ be a unital associative $\mathbb{K}$-algebra as before.

\begin{definition}
Let $M$ be a left $B$-module. 
\begin{itemize}
	\item A \textit{connection} is a $\mathbb{K}$-linear map $\nabla\colon M\to \Omega^1B\otimes_BM$ satisfying the Leibniz rule:
	\[
		\nabla(b\cdot m) = db\otimes m + b\cdot \nabla(m)\quad \textrm{for } b\in B \textrm{ and } m\in M\,.
	\]
	This extends to a degree $1$ derivation on $\Omega^\bullet B\otimes_BM$ by
	\[
		\nabla(\omega\otimes m) = d\omega\otimes m + (-1)^p \omega\cdot\nabla(m)\quad \textrm{for } \omega\in \Omega^pA \textrm{ and } m\in M\,.
	\]
	\item The \textit{curvature} of a connection $\nabla$ is defined by $R=\nabla^2\colon M\to \Omega^2B\otimes_BM$, which is a $B$-module map.
	\item A connection $\nabla$ is \textit{flat} if the curvature $R$ vanishes identically.
\end{itemize}
\end{definition}

\begin{remark}\label{rem:connproj}
Since $M$ is $\mathbb{K}$-free, $M$ admits a ($\mathbb{K}$-linear) connection if and only if $M$ is $B$-projective (see, for example, Corollary 8.2 of \cite{cq}). The existence of flat connections is more subtle; one sufficient condition is that $M$ is $B$-free. In fact, a connection on a free module is uniquely specified by its values on a free basis. One necessary condition, on the other hand, is the vanishing of the Chern classes; for the definition, see \cite{Ginzburg:aa}.\\
\end{remark}

\begin{definition}
Let $\nabla$ be a connection on a $B$-module $M$. Its \textit{push-out} by a $\mathbb{K}$-algebra homomorphism $\psi:B\to C$ is the connection $\psi_*\nabla$ on $C\otimes_BM$ defined as follows: first of all, we have the induced map $\psi\otimes\id\colon\Omega^1B\otimes_B M\to \Omega^1C\otimes_B M$.
Then, we define $\psi_*\nabla$ by
\begin{align*}
	\psi_*\nabla \colon C\otimes_B M &\to \Omega^1C\otimes_BM\\
	c\otimes m &\mapsto dc\otimes m + c\psi(\nabla(m))
\end{align*}
with the target naturally identified with $\Omega^1C\otimes_{C}(C\otimes_B M)$\,.
\end{definition}

Next, we recall the trace of a module endomorphism over a non-commutative algebra. Firstly, we have the evaluation map
\begin{align*}
	\mathrm{|ev|}&\colon M^*\otimes_B M \to |B|\colon \theta\otimes m\mapsto |\theta(m)|\,.
\end{align*}
If $M$ is a dualisable $B$-module, we have the isomorphism $\End_B(M) \cong M^*\otimes_BM$ as seen in Proposition \ref{prop:dualisable}, and the composite
\[
	\Tr \colon \End_B(M)\xrightarrow{\sim} M^*\otimes_B M \xrightarrow{|\mathrm{ev}|} |B|\,,
\]
is known as the \textit{Hattori--Stallings trace} \cite{httrace,STALLINGS1965170}.
\begin{remark}\label{rem:trace}
Explicitly, once we take maps $(e,e^*)$ by  Proposition \ref{prop:dualisable}, the trace of $\psi\in\End_B(M)$ is given by
\[
	\Tr(\psi) = \sum_i e^*_i(\psi(e_i))\,.
\]
\end{remark}

\begin{proposition}\label{prop:HTtraceLie}
Let $M$ and $N$ be dualisable $B$-modules. Then, for $B$-module homomorphisms $f:M\to N$ and $g:N\to M$, we have $\Tr(f\circ g) = \Tr(g\circ f)$. 
\end{proposition}
\noindent Proof. For $f=\theta\otimes n\in M^*\otimes_B N$ and $g=\varphi\otimes m\in N^*\otimes_B M$, we have
\begin{align*}
	f\circ g = \varphi\otimes\theta(m)n\;\textrm{ and }\; g\circ f = \theta\otimes\varphi(n)m\,,
\end{align*}
and therefore
\begin{align*}
	\Tr(f\circ g) &= |\theta(m)\varphi(n)| = |\varphi(n)\theta(m)| = \Tr(g\circ f).
\end{align*}
This completes the proof.\qed\\

We need a suitable action of derivations on $M$ to define the associated divergence. 

\begin{definition}\label{def:deract}
A \textit{derivation action} on a left $B$-module $M$ is a triple $(\mathfrak{d},\varphi,\rho)$, where $\mathfrak{d}$ is a Lie algebra, and $\varphi\colon \mathfrak d\to \Der_\mathbb{K}(B)$ and $\rho\colon\mathfrak{d}\to \End_\mathbb{K}(M)$ are Lie algebra homomorphisms satisfying the compatibility condition:
\begin{align}
	\rho(f)(bm) = \varphi(f)(b)\cdot m + b\cdot \rho(f)(m) \quad \textrm{for } f\in \mathfrak{d}, b\in B \textrm{ and } m\in M.\label{eq:compat}
\end{align}
\end{definition}
\noindent Given a homomorphism $\varphi$, $|B|$ is naturally a $\mathfrak d$-module by the composite
\[
	\mathfrak d\xrightarrow{\varphi} \Der_\mathbb{K}(B) \to \End_\mathbb{K}(|B|)
\]
of Lie algebra homomorphisms.

\begin{definition}\label{def:conn}
	Let $\nabla$ be a connection on a dualisable $M$. The \textit{divergence map} associated with $\nabla$ and a derivation action $(\mathfrak{d},\varphi,\rho)$ is defined by 
	\[
		\Div^{(\nabla,\varphi,\rho)}\colon\mathfrak d\to |B|\colon f\mapsto \Tr(\rho(f) - (i_{\varphi(f)}\otimes\id_M)\circ\nabla)\,.
	\]
\end{definition}
\noindent Note that $\rho(f) - (i_{\varphi(f)}\otimes\id_M)\circ\nabla$ is a $B$-module homomorphism by the compatibility condition.\\

Now assume that $A$ is formally smooth until the end of this section, so that $\Omega^1A$ is dualisable and
\[
	\DDer_\mathbb{K}(A)\otimes_{A^\mathrm{e}}\Omega^1A\cong \End_{A^\mathrm{e}}(\Omega^1A)
\]
holds. The trace map now takes the form $\Tr\colon \End_{A^\mathrm{e}}(\Omega^1A)\to |A^\mathrm{e}|$. Setting $B = A^\mathrm{e}$ and $M = \Omega^1A$, we have the default derivation action given by $\mathfrak d = \Der_\mathbb{K}(A)$,
\begin{align}
	\begin{aligned}\label{eq:defaultaction}
	\varphi_\mathrm{assoc}\colon \Der_\mathbb{K}(A)& \to \Der_\mathbb{K}(A^\mathrm{e})\hspace{9.5pt}\colon \;f \mapsto  \tilde f := f\otimes \id_{A^\mathrm{op}} + \id_A\otimes f,\\
	\rho_\mathrm{assoc}\colon  \Der_\mathbb{K}(A) &\to  \End_\mathbb{K}(\Omega^1A)\colon \;f\mapsto L_f\,.
	\end{aligned}
\end{align}
The compatibility condition is automatically satisfied as $L_f$ is a degree 0 derivation on $\Omega^\bullet A$. Thus, we obtain the divergence
\[
	\Div^\nabla\colon \Der_\mathbb{K}(A)\to |A^\mathrm{e}|\colon f \mapsto \Tr(L_f - (i_{\tilde f}\otimes\id_{\Omega^1A})\circ\nabla)
\]
associated with a connection $\nabla\colon\Omega^1A\to \Omega^1A^\mathrm{e}\otimes_{A^\mathrm{e}}\Omega^1A$ and the default derivation action. In this case, we have an interpretation of $\Div^\nabla$ in terms of the horizontal lift:
\begin{deflem}\label{deflem:horizontal}
For a connection $\nabla\colon\Omega^1A\to \Omega^1A^\mathrm{e}\otimes_{A^\mathrm{e}}\Omega^1A$, there is a unique $\mathbb{K}$-linear map
\[
	(\,\cdot\,)^\mathrm{H}\colon\Der_\mathbb{K}(A)\to  \Der_\mathbb{K}(\Omega^\bullet A)^{(0)} \colon f\mapsto f^\mathrm{H}
\]
satisfying the following properties: for $f\in \Der_\mathbb{K}(A)$,
\begin{enumerate}[(1)]
	\item $f^\mathrm{H}(a) = f(a)$ in $A = \Omega^0A$ for $a\in A$; and
	\item $f^\mathrm{H}(\alpha) =(i_{\tilde f}\otimes\id_{\Omega^1A})\circ\nabla(\alpha)$ in $A^\mathrm{e}\otimes_{A^\mathrm{e}}\Omega^1A\cong \Omega^1A$ for $\alpha\in \Omega^1A$.
\end{enumerate}
The map $(\,\cdot\,)^\mathrm{H}$ is called the horizontal lift by $\nabla$.
\end{deflem}
\noindent Proof. The uniqueness follows from the fact that $\Omega^\bullet A$ is generated as a $\mathbb{K}$-algebra by $A$ and $\Omega^1A$. Then, as $\Omega^\bullet A$ is the tensor algebra of $\Omega^1A$ over $A$, it suffices to show that (1) and (2) are compatible to see that $f^\mathrm{H}$ is well-defined.  We compute, for $x,y\in A$ and $\alpha\in\Omega^1A$,
\begin{align*}
	f^\mathrm{H}(x\alpha y) &= (i_{\tilde f}\otimes\id_{\Omega^1A})\circ\nabla(x\bar y\cdot\alpha) &\textrm{by the definition of the horizontal lift,}\\
	&= (i_{\tilde f}\otimes\id_{\Omega^1A})(d(x\bar y)\otimes\alpha + x\bar y\cdot\nabla(\alpha)) &\textrm{since }\nabla\textrm{ is a connection,}\\
	&= \tilde f(x\bar y)\cdot \alpha +  x\bar y\cdot f^\mathrm{H}(\alpha)&\textrm{by the definition of the contraction map,}\\
	&= f(x)\alpha y + xf^\mathrm{H}(\alpha)y + x\alpha f(y)&\textrm{by the definition of }\tilde f\,.
\end{align*}
This completes the proof.\qed\\

\noindent For the original definition of the horizontal lift in Riemannian geometry, see \cite{horizontal}.\\

There is a simple criterion for $\Div^\nabla$ to be a Lie algebra $1$-cocycle with coefficients in $|A^\mathrm{e}|$.

\begin{proposition}\label{prop:flatcocyc}
If $\nabla$ is a flat connection, then $\Div^\nabla \!= \Tr(L - (\,\cdot\,)^\mathrm{H})$ is a Lie algebra $1$-cocycle.
\end{proposition}

We need some preparation for the proof. First, define the action of $\Der_\mathbb{K}(A)$ on $\DDer_\mathbb{K}(A)$ by
\[f
	\cdot \theta = [f,\theta] := (f\otimes\id + \id\otimes f)\circ\theta - \theta\circ f\quad \textrm{for } f\in \Der_\mathbb{K}(A)\textrm{ and } \theta\in\DDer_\mathbb{K}(A)
\]
with the commutator taken in the space $\Der_\mathbb{K}(T(A))$, where $T(A)$ denotes the tensor algebra generated by $A$ over $\mathbb{K}$. Next, define the action of $\Der_\mathbb{K}(A)$ on $\End_{A^\mathrm{e}}(\Omega^1A)$ by
\[
	f\cdot \psi = [L_f,\psi]\quad\textrm{ for } f\in\Der_\mathbb{K}(A)\textrm{ and } \psi\in\End_{A^\mathrm{e}}(\Omega^1A)\,.
\]
With these actions, the isomorphism of $\mathbb{K}$-modules $\DDer_\mathbb{K}(A)\otimes_{A^\mathrm{e}}\Omega^1A\cong \End_{A^\mathrm{e}}(\Omega^1A)$ becomes one of $\Der_\mathbb{K}(A)$-modules.

\begin{lemma}\label{lem:HTtraceDer}
The trace map $\Tr$ is a $\Der_\mathbb{K}(A)$-module homomorphism.
\end{lemma}
\noindent Proof. For $\theta\otimes da \in\DDer_\mathbb{K}(A)\otimes_{A^\mathrm{e}}\Omega^1A\cong \End_{A^\mathrm{e}}(\Omega^1A)$ and $f\in\Der_\mathbb{K}(A)$, we have
\begin{align*}
	f\cdot(\theta\otimes da)& = [f,\theta]\otimes da + \theta\otimes L_f(da),
\end{align*}
by the definition of the action, so that
\begin{align*}
	\Tr(f\cdot(\theta\otimes da)) &= |[f,\theta](a) + \theta(f(a))| = |f(\theta(a))| = f\cdot\Tr(\theta\otimes da).
\end{align*}
The $\Der_\mathbb{K}(A)$-equivariance of the trace follows. \qed\\

\begin{lemma}\label{lem:LiehomR}
Let $R$ be the curvature of $\nabla$ and $f,g\in \Der_\mathbb{K}(A)$. Then we have $[f,g]^\mathrm{H} = [f^\mathrm{H},g^\mathrm{H}] + i_{\tilde f}i_{\tilde g}R$ in $\End_\mathbb{K}(\Omega^1A)$.
\end{lemma}
\noindent Proof. For $\alpha\in\Omega^1A$, put $\nabla\alpha = \sum\omega\otimes\beta$ for some $\omega\in\Omega^1A^\mathrm{e}$ and $\beta\in\Omega^1A$. Dropping the summation symbol and denoting $i_{\tilde f}\otimes\id_{\Omega^1A}$ simply by $i_{\tilde f}$, we have
\begin{align*}
	[f^\mathrm{H},g^\mathrm{H}]\alpha &= f^\mathrm{H}(i_{\tilde g}\nabla\alpha)- g^\mathrm{H}(i_{\tilde f}\nabla\alpha)\\
	&= f^\mathrm{H}(i_{\tilde g}(\omega\otimes\beta))- g^\mathrm{H}(i_{\tilde f}(\omega\otimes\beta))\\
	&= f^\mathrm{H}(i_{\tilde g}\omega\cdot\beta)- g^\mathrm{H}(i_{\tilde f}\omega\cdot\beta)\,.
\end{align*}
Since $f^\mathrm{H}$ and $g^\mathrm{H}$ are derivations by Definition-Lemma \ref{deflem:horizontal}, this is equal to
\begin{align*}
	&\, \tilde f(i_{\tilde g}\omega)\cdot\beta +i_{\tilde g}\omega\cdot f^\mathrm{H}(\beta) - \tilde g(i_{\tilde f}\omega)\cdot\beta - i_{\tilde f}\omega\cdot g^\mathrm{H}(\beta)\\
	=&\, \tilde f(i_{\tilde g}\omega)\cdot\beta +i_{\tilde g}\omega\cdot i_{\tilde f}\nabla\beta - \tilde g(i_{\tilde f}\omega)\cdot\beta - i_{\tilde f}\omega\cdot i_{\tilde g}\nabla\beta
\end{align*}
together with the definition of the horizontal lift. On the other hand, for the curvature, we have 
\begin{align*}
	R(\alpha) &= \nabla(\omega\otimes\beta) = d\omega\otimes\beta - \omega\nabla\beta, \textrm{ and}\\
	i_{\tilde f}i_{\tilde g}R(\alpha) &= i_{\tilde f}i_{\tilde g}d\omega\cdot\beta - i_{\tilde g}\omega \cdot i_{\tilde f}\nabla\beta + i_{\tilde f}\omega \cdot i_{\tilde g}\nabla\beta.
\end{align*}
Therefore, the right-hand side of the claim is equal to 
\begin{align*}
	\left([f^\mathrm{H},g^\mathrm{H}] + i_{\tilde f}i_{\tilde g}R\right)(\alpha) &= (\tilde f(i_{\tilde g}\omega)  - \tilde g(i_{\tilde f}\omega) + i_{\tilde f}i_{\tilde g}d\omega)\cdot\beta\\
	&= [i_{\tilde f},L_{\tilde g}]\omega\cdot\beta = i_{[\tilde f,\tilde g]}\omega\cdot\beta\\
	&=  i_{[\tilde f,\tilde g]}\nabla\alpha = [f,g]^\mathrm{H}(\alpha).
\end{align*}
This completes the proof.\qed\\

\noindent\textbf{Proof of Proposition \ref{prop:flatcocyc}.} Consider the following split sequence of Lie algebras:
\[\begin{tikzcd}[cramped]
	0 \arrow[r] & \Der_A(\Omega^\bullet A)^{(0)} \arrow[r, "\mathrm{incl}"] & \Der_\mathbb{K}(\Omega^\bullet A)^{(0)} \arrow[r,shift left,"\mathrm{res}_A"] &  \Der_\mathbb{K}(A) \arrow[l,shift left,"L"]\arrow[r] & 0\,,
\end{tikzcd}\]
where $\mathrm{res}_A$ is the restriction to $A$. Indeed, if two elements in $ \Der_\mathbb{K}(\Omega^\bullet A)^{(0)}$ restrict to the same derivation on $A$, their difference is an $A$-linear derivation, and therefore contained in the image of $\mathrm{incl}$.

Since any $A$-linear derivation on $\Omega^\bullet A$ is uniquely determined by its restriction to $\Omega^1A$, the space $\Der_A(\Omega^\bullet A)^{(0)}$ is isomorphic to $\End_{A^\mathrm{e}}(\Omega^1A)$ as a Lie algebra. Then we define
\[
	c\colon \Der_\mathbb{K}(A) \to \End_{A^\mathrm{e}}(\Omega^1A): f \mapsto f^\mathrm{H}-L_f
\]
using this isomorphism. By Lemma \ref{lem:LiehomR} and the flatness of $\nabla$, $(\,\cdot\,)^\mathrm{H}$ is a Lie algebra homomorphism. Then $c$ satisfies the non-abelian $1$-cocycle condition:
\begin{align}\label{eq:na1cocyc}
	c([f,g]) = f\cdot c(g) - g\cdot c(f) + [c(f),c(g)]\,.
\end{align}
To check this, we compute
\begin{align*}
	c([f,g]) &= [f,g]^\mathrm{H} - L_{[f,g]}\\
	&= [f^\mathrm{H},g^\mathrm{H}] - [L_f,L_g]\\
	&= [c(f),g^\mathrm{H}] + [L_f,g^\mathrm{H}]- [L_f,L_g]\\
	&= [c(f),c(g)] + [c(f),L_g] + [L_f,g^\mathrm{H} - L_g]\\
	&= [c(f),c(g)] - [L_g,c(f)] + [L_f,c(g)]\\
	&= [c(f),c(g)] - g\cdot c(f) + f\cdot c(g)\,.
\end{align*}
Finally, applying the trace to (\ref{eq:na1cocyc}) and using Proposition \ref{prop:HTtraceLie} and Lemma \ref{lem:HTtraceDer}, $\Div^\nabla=-\Tr\circ\,c$ is indeed a Lie algebra $1$-cocycle.\qed

\begin{remark}
Proposition \ref{prop:flatcocyc} remains true for the general case $\Div^{(\nabla,\varphi,\rho)}$ in Definition \ref{def:conn}; the proof is almost identical, just less conceptual.
\end{remark}

Similar formulations of a connection and the associated divergence map in differential geometry can be seen in the context of Lie algebroid theory. For example, a connection is defined using an analogous exact sequence to the above in Definition 3.8 of \cite{Lazzarini:2010aa}, and the corresponding divergence formula can be found in Proposition 3.11 of \cite{Xu1999}.

An operad-theoretic formulation of the divergence map, on the other hand, has been studied in \cite{Powell:2021aa}, where they defined the standard divergence map on a free $\mathscr{O}$-algebra for a reduced operad $\mathscr{O}$. The standard divergences are also discussed in the next two sections, only in the case of operads controlling associative algebras and Lie algebras. As we will see below, they correspond to those associated with the canonical flat connections naturally defined from free-generating systems in our formulation.\\

\section{Examples of Flat Connections}\label{sec:flatconns}

In this short section, we look into two examples of flat connections.

Firstly, we consider the case $A = \mathbb{K}\langle z_1,\dotsc,z_r\rangle$. We define the standard connection $\nabla\!_z$  on $\Omega^1A$ associated with the free-generating system $z = (z_i)_{1\leq i\leq r}$  by $\nabla\!_z(dz_i) = 0$ for all $i$. Then, this connection is obviously flat. Now, since $(dz_i)_{1\leq i\leq r}$ is a $A^\mathrm{e}$-free basis of $\Omega^1A$, we can take its dual basis $(\partial_i)_{1\leq i\leq r}$, which comprises of double derivations. More precisely, they are given by the formula
\[
	\partial_j(z_i) = \delta_{ij}\cdot 1\otimes 1\quad\textrm{ for } 1\leq i,j\leq r.
\]
For $f \in\Der_\mathbb{K}(A)$, we compute
\begin{align*}
	(L_f-f^\mathrm{H})(dz_i) &= df(z_i) - i_{\tilde f}\nabla\!_z(dz_i) = \sum_j\partial_j(f(z_i))\cdot dz_j,\textrm{ and}\\
	\Div^{\nabla\!_z}(f) &= \Tr(L_f-f^\mathrm{H}) = \sum_i|\partial_i(f(z_i))|\,,
\end{align*}
which recovers the standard (double) divergence map.\\

Next, we investigate the case $A = \mathbb{K}\pi$. Recall that the divergence $\Div^\mathcal{C}$ associated with a free-generating system $\mathcal{C} = (a_i,b_i,\zeta_j)_{1\leq i\leq g,1\leq j\leq n}$, introduced in (\ref{eq:divc}), is given by
\[
	\Div^\mathcal{C}(f) = \sum_{c\in\mathcal{C}} |\partial_c(f(c)) - 1\otimes c^{-1}f(c)|\in|A|\otimes|A|.
\]
Here $\partial_c$ is given by $\partial_c(c) = \delta_{c,c'} \cdot 1\otimes 1$. Since $(dc\,c^{-1})_{c\in\mathcal{C}}$ is a $A^\mathrm{e}$-free basis of $\Omega^1A$ as we proved in Proposition \ref{prop:freegroupsmooth}, we define the connection $\nabla_\mathcal{C}$ associated with $\mathcal{C}$ by the formula $\nabla_\mathcal{C}(dc\,c^{-1}) = 0$ for all $c\in\mathcal{C}$, which is also flat. The equality $\nabla_\mathcal{C}(dc\,c^{-1}) = 0$ is equivalent to
\begin{align*}
	0 = \bar c^{-1}\nabla_\mathcal{C}(dc) + d\bar c^{-1}\otimes dc\,,\:\;\textrm{ or }\;\;\nabla_\mathcal{C}(dc) = d\bar c  \,\bar c^{-1}\otimes dc.
\end{align*}
Here we used $d(c^{-1}) = -c^{-1}dc\,c^{-1}$.

Now, we compute the associated divergence. We know that $(dc)_{c\in\mathcal{C}}$ is also a free $A^\mathrm{e}$-basis of $\Omega^1A$ since $c$ is invertible in $A$. Then, we have
\begin{align*}
	(L_f-f^\mathrm{H})(dc) &= df(c) - i_{\tilde f}\nabla_\mathcal{C}(dc)\\
	&= \Big(\sum_{c'\in\mathcal{C}} \partial_{c'}(f(c))\cdot dc' \Big) - (f(\bar c)\bar c^{-1})\cdot dc,\textrm{ and}\\
	\Div^{\nabla_\mathcal{C}}(f) &= \sum_{c\in\mathcal{C}}|\partial_c(f(c)) - f(\bar c)\bar c^{-1}|,
\end{align*}
which coincides with $\Div^\mathcal{C}(f)$ using the identification $|A^\mathrm{e}|\cong|A|\otimes|A|:|x\bar y| \mapsto |x|\otimes |y|$.\\

\section{Connections on Hopf Algebras and Geometry over the Lie Operad}\label{sec:Lie}

The two connections in the previous section are induced by connections involving geometry over the Lie operad, as shown below. However, we will not use any operad explicitly. We first discuss a correspondence of connections over a Hopf algebra and then apply it to the enveloping algebra of a Lie algebra. At the end of this section, we briefly discuss the relation between the divergence map and the Satoh trace. This section may be skipped since this material is not used in the following sections.\\

Let $(A,\mu,\Delta,\eta,\varepsilon,S)$ be a Hopf algebra and denote by
\[
	\Der_\mathrm{Hopf}(A) = \{f\in\Der_\mathbb{K}(A)\colon (f\otimes\id + \id\otimes f)\circ\Delta = \Delta\circ f\}
\]
the space of Hopf algebra derivations. We have two algebra homomorphisms:
\[
	\tilde\Delta := (\id\otimes\, S)\circ\Delta\colon A\to A^\mathrm{e}\;\textrm{ and }\;\id\otimes\, \varepsilon \colon A^\mathrm{e}\to A\,.
\]
Recall the functor $\Phi_A: A\textrm{-}\textbf{Mod}\to A^\mathrm{e}\textrm{-}\textbf{Mod}:M\mapsto M\otimes A$ in Definition \ref{def:functorphi}. 

\begin{deflem}
Define the natural transformation $j\colon\id_{A\textrm{\rm-}\mathbf{Mod}}\to\Phi_A$ by, for an $A$-module $M$,
\begin{align*}
	j_M\colon M\to\Phi_A(M): m\mapsto m\otimes 1.
\end{align*}
Then, each $j_M$ is a module homomorphism over $\tilde\Delta$.
\end{deflem}
\noindent Proof. It is clear that $j$ defines a natural transformation. Next, we check that $j_M$ is a module homomorphism over $\tilde\Delta$. For $m\in M$ and $x\in A$, we have
\begin{align*}
	\tilde\Delta(x)\cdot j_M(m) &= x^{(1)}\cdot (m\otimes 1)\cdot S(x^{(2)})\\
	&= x^{(1)}m\otimes x^{(2)}S(x^{(3)})
\end{align*}
by the definition of the $A^\mathrm{e}$-action on $\Phi_A(M)$. Since we have
\[
	x^{(1)}\otimes x^{(2)}S(x^{(3)}) = x^{(1)}\otimes \varepsilon(x^{(2)}) = x\otimes 1
\]
for a Hopf algebra by the defining property of the antipode, we obtain
\[
	\tilde\Delta(x)\cdot j_M(m) = xm\otimes 1 = j_M(xm).
\]
This completes the proof. \qed\\

Next, we show that the functor $\Phi_A$ extends to the space of connections. Recall the notation $x\bar y\in A^\mathrm{e}$ where $x,y\in A$. With this notation, $dx\,\bar y$ and $\bar y dx$ are elements of $\Omega^1A^\mathrm{e}$, but they are not necessarily equal.

\begin{deflem}
Let $M$ be an $A$-module and $\nabla\colon M\to\Omega^1A\otimes M$ be a connection. Then, there is a unique connection $\Phi_A(\nabla)$ on $\Phi_A(M)$ that makes the following diagram commute:
\[\begin{tikzcd}[ampersand replacement = \&]
	M\arrow[r, "\nabla"] \arrow[d,"j_M"]\& \Omega^1A\otimes_A M\arrow[d,"\tilde\Delta\otimes j_M"]\\
	\Phi_A(M)\arrow[r, "\Phi_A(\nabla)"] \& \Omega^1A^\mathrm{e}\otimes_{A^\mathrm{e}}\Phi_A(M)
\end{tikzcd}\]
Furthermore, $\Phi_A$ preserves flat connections.
\end{deflem}
\noindent Proof. We define $\Phi_A(\nabla)$ by the formula
\[
	\Phi_A(\nabla)(m\otimes a) = d\bar a\otimes (m\otimes 1) + \bar a\cdot (\tilde\Delta\otimes j_M)(\nabla(m))\,.
\]
Then, the square above is clearly commutative. We check that this defines a connection. For $m\in M$ and $x,y,a\in A$, we have
\begin{align*}
	\Phi_A(\nabla)(x\cdot (m\otimes a)\cdot y) &= \Phi_A(\nabla)(x^{(1)}m\otimes x^{(2)}ay)\\
	&= d(\overline{x^{(2)}ay})\otimes (x^{(1)}m\otimes 1) + \overline{x^{(2)}ay} \cdot (\tilde\Delta\otimes j_M)(\nabla(x^{(1)}m))\,.
\end{align*}
The first term is equal to
\begin{align*}
	 d(\overline{x^{(2)}ay})\otimes (x^{(1)}m\otimes 1) &= d(\overline{x^{(3)}ay}) x^{(1)}\otimes (m\otimes S(x^{(2)}))\\
	 &= d(\overline{x^{(3)}ay}) x^{(1)}\overline{S(x^{(2)})}\otimes (m\otimes 1)\,,
\end{align*}
while the second term is equal to
\begin{align*}
	\overline{x^{(2)}ay} \cdot (\tilde\Delta\otimes j_M)(\nabla(x^{(1)}m)) &= \overline{x^{(2)}ay} \cdot (\tilde\Delta\otimes j_M)( dx^{(1)} \otimes m + x^{(1)}\cdot\nabla(m) )\\
	&= \overline{x^{(2)}ay} \cdot \big[ d\tilde\Delta(x^{(1)}) \otimes (m\otimes 1) + \tilde\Delta(x^{(1)})\cdot (\tilde\Delta\otimes j_M)(\nabla(m)) \big]\\
	&= \overline{x^{(3)}ay} \cdot \big[ d(x^{(1)}\overline{S(x^{(2)})}) \otimes (m\otimes 1) + x^{(1)}\overline{S(x^{(2)})}\cdot (\tilde\Delta\otimes j_M)(\nabla(m)) \big]\,.
\end{align*}
The sum of the coefficients of $m\otimes 1$ in the above two is now equal to
\begin{align*}
	d(\overline{x^{(3)}ay}) x^{(1)}\overline{S(x^{(2)})} +  \overline{x^{(3)}ay} \,d(x^{(1)}\overline{S(x^{(2)})}) &= d(\overline{x^{(3)}ay}\,x^{(1)}\overline{S(x^{(2)})})\\
	&= d(x^{(1)}\overline{S(x^{(2)})x^{(3)}ay}) = d(x\overline{ay})\,.
\end{align*}
In addition, the coefficient of $(\tilde\Delta\otimes j_M)(\nabla(m))$ is equal to
\[
	\overline{x^{(3)}ay}x^{(1)}\,\overline{S(x^{(2)})} = x^{(1)}\overline{S(x^{(2)})x^{(3)}ay} = x\overline{ay}\,.
\]
Therefore, we obtain
\begin{align*}
	\Phi_A(\nabla)(x\cdot (m\otimes a)\cdot y) &= d(x\overline{ay})\otimes (m\otimes 1) + x\overline{ay}\cdot (\tilde\Delta\otimes j_M)(\nabla(m))\\
	&= d(x\bar y\bar a)\otimes (m\otimes 1) + x\bar y\bar a\cdot (\tilde\Delta\otimes j_M)(\nabla(m))\\
	&= (d(x\bar y)\bar a + x\bar yd\bar a)\otimes (m\otimes 1) + x\bar y\bar a\cdot (\tilde\Delta\otimes j_M)(\nabla(m))\\
	&= d(x\bar y)\otimes (m\otimes a) + x\bar y\cdot \Phi_A(\nabla)(m\otimes a)\,,
\end{align*}
which shows that $\Phi_A(\nabla)$ is a connection. The uniqueness follows from the fact that the image of $j_M$ generates $\Phi_A(M)$ as an $A^\mathrm{e}$-module.

Finally, if $\nabla$ is flat, the composite
\[
	(\tilde\Delta\otimes j_M)\circ\nabla\circ\nabla = \Phi_A(\nabla)\circ\Phi_A(\nabla)\circ j_M
\]
is zero. This forces $\Phi_A(\nabla)^2$ to be zero since, again, the image of $j_M$ generates $\Phi_A(M)$ as an $A^\mathrm{e}$-module.\qed\\

We can extend a derivation action along $\Phi_A$ under a suitable assumption.

\begin{deflem}
Let $(\mathfrak d,\varphi,\rho)$ be a derivation action on an $A$-module $M$. If the image of $\varphi\colon \mathfrak d\to\Der_\mathbb{K}(A)$ is contained in $\Der_\mathrm{Hopf}(A)$, we have the induced derivation action $(\mathfrak d,\Phi_A(\varphi),\Phi_A(\rho))$ on $\Phi_A(M)$ given by
\begin{align*}
	\Phi_A(\varphi)&\colon \mathfrak d\xrightarrow{\varphi}\Der_\mathbb{K}(A)\xrightarrow{\varphi_\mathrm{assoc}} \Der_\mathbb{K}(A^\mathrm{e})\,,\textrm{ and}\\
	\Phi_A(\rho)&\colon \mathfrak d\to\End_\mathbb{K}(\Phi_A(M)): f \mapsto (m\otimes a \mapsto \rho(f)(m)\otimes a + m\otimes \varphi(f)(a))
\end{align*}
\end{deflem}
\noindent Proof. We check the compatibility condition (\ref{eq:compat}). For $f\in\mathfrak d$, $m\in M$ and $x,y,a\in A$, the left-hand side of (\ref{eq:compat}) is equal to
\begin{align*}
	\Phi_A(\rho)&(f)(x\cdot(m\otimes a)\cdot y)\\
	&= \Phi_A(\rho)(f)(x^{(1)}m\otimes x^{(2)}a y)\\
	&= \rho(f)(x^{(1)}m)\otimes x^{(2)}a y + x^{(1)}m\otimes \varphi(f)(x^{(2)}a y)\\
	&= \varphi(f)(x^{(1)})m\otimes x^{(2)}a y + x^{(1)}\rho(f)(m)\otimes x^{(2)}a y \\
	&\qquad + x^{(1)}m\otimes \varphi(f)(x^{(2)})ay + x^{(1)}m\otimes x^{(2)}\varphi(f)(a)y + x^{(1)}m\otimes x^{(2)}a \varphi(f)(y)\,.
\end{align*}
On the other hand, the right-hand side of (\ref{eq:compat}) is equal to,
\begin{align*}
	\Phi_A(\varphi)&(f)(x\bar y)\cdot (m\otimes a)+ (x \bar y)\cdot \Phi_A(\rho)(f)(m\otimes a)\\
	&= \varphi(f)(x)\cdot (m\otimes a)\cdot y + x\cdot \Phi_A(\rho)(f)(m\otimes a)\cdot y + x\cdot (m\otimes a)\cdot \varphi(f)(y)\\
	&=  \varphi(f)(x)^{(1)}m\otimes \varphi(f)(x)^{(2)}ay + x^{(1)}\rho(f)m\otimes x^{(2)}ay + x^{(1)}m\otimes x^{(2)}\varphi(f)(a)y+ x^{(1)}m\otimes x^{(2)}a \varphi(f)(y)\,.
\end{align*}
Therefore, they are equal if the equality
\[
	\varphi(f)(x^{(1)})\otimes x^{(2)} + x^{(1)}\otimes \varphi(f)(x^{(2)}) =  \varphi(f)(x)^{(1)}\otimes \varphi(f)(x)^{(2)}
\]
holds. This equality follows from $(\varphi(f)\otimes\id + \id\otimes\,\varphi(f)) \circ\Delta = \Delta\circ\varphi(f)$, which, in turn, follows from the assumption that $\varphi(f)$ is a Hopf algebra derivation.\qed\\

The correspondence between connections is reflected in the associated divergences.
\begin{theorem}\label{thm:hopfconn}
Let $M$ be a dualisable $A$-module, $\nabla\colon M\to\Omega^1A\otimes_A M$ a connection, and $(\mathfrak d,\varphi,\rho)$ a derivation action on $M$ with $\Imag\varphi\subset\Der_\mathrm{Hopf}(A)$. Then, we have
\[
	\Div^{(\Phi_A(\nabla),\Phi_A(\varphi),\Phi_A(\rho))} = \tilde\Delta\circ \Div^{(\nabla,\varphi,\rho)}
\]
as maps $\mathfrak d \to |A^\mathrm{e}|$.
\end{theorem}

\begin{lemma}
Let $M$ be a dualisable $A$-module, $e\colon A^{\oplus r}\to M$ and $e^*\colon M\to A^{\oplus r}$ the maps in Proposition \ref{prop:dualisable}, and $F\colon \Phi_A(A)\to A^\mathrm{e}$ the isomorphism in (\ref{eq:bimodisom}).
\begin{enumerate}[(1)]
	\item The $A^\mathrm{e}$-module maps $\Phi_A(e)\circ (F^{-1})^{\oplus r}\colon {A^\mathrm{e}}^{\oplus r}\to \Phi_A(M)$ and $F^{\oplus r}\circ \Phi_A(e^*)\colon \Phi_A(M)\to {A^\mathrm{e}}^{\oplus r}$ satisfy
	\[
		 \Phi_A(e)\circ (F^{-1})^{\oplus r}\circ F^{\oplus r}\circ \Phi_A(e^*) = \id_{\Phi_A(M)}.
	\]
	\item $\Phi_A(M)$ is a dualisable $A^\mathrm{e}$-module.
\end{enumerate}
\end{lemma}
\noindent Proof. (1) It is clear from the condition $e\circ e^* = \id_M$. (2) follows from (1) together with Proposition \ref{prop:dualisable}.\qed\\

\noindent \textbf{Proof of Theorem \ref{thm:hopfconn}.} By (2) of the lemma above, talking about the trace map (and hence divergence maps) on $\Phi_A(M)$ makes sense. Define $e^i$ and $e^*_i$ as in the proof of Proposition \ref{prop:dualisable},
\begin{align*}
	\varepsilon^i := \Phi_A(e)\circ (F^{-1})^{\oplus r}(0,\dotsc,1\otimes 1,\dotsc,0) = e^i\otimes 1 = j_M(e^i)\,,
\end{align*}
and $\varepsilon^*_i$ to be the $i$-th component of $F^{\oplus r}\circ \Phi_A(e^*)$. Then, we have $\varepsilon^*_i\circ j_M = F\circ j_A\circ e^*_i$ since $j$ is a natural transformation. For $f\in\mathfrak d$, we compute,
\begin{align*}
	&\Div^{(\Phi_A(\nabla),\Phi_A(\varphi),\Phi_A(\rho))}(f)&\\
	&= \sum_i \Big| \varepsilon^*_i\big[\Phi_A(\rho)(f)(\varepsilon^i) - (i_{\Phi_A(\varphi)(f)}\otimes\id)\circ\Phi_A(\nabla)(\varepsilon^i)\big] \Big|& \textrm{by Remark \ref{rem:trace},} \\
	&= \sum_i \Big| \varepsilon^*_i\big[\Phi_A(\rho)(f)(j_M(e^i)) - (i_{\Phi_A(\varphi)(f)}\otimes\id)\circ\Phi_A(\nabla)(j_M(e^i))\big] \Big|& \textrm{by the definition of } \varepsilon^i,\\
	&= \sum_i \Big| \varepsilon^*_i\big[j_M(\rho(f)(e^i)) - (i_{\Phi_A(\varphi)(f)}\otimes\id)\circ(\tilde\Delta\otimes j_M)(\nabla(e^i))\big] \Big|& \textrm{by the definition of } \Phi_A(\rho) \textrm{ and } \Phi_A(\nabla),\\
	&= \sum_i \Big| \varepsilon^*_i\big[j_M(\rho(f)(e^i)) - j_M\circ (i_{\varphi(f)}\otimes\id)(\nabla(e^i))\big] \Big|& \textrm{since }\varphi(f)\textrm{ is a Hopf derivation,}\\
	&= \sum_i \Big| F\circ j_A\circ e^*_i\big[\rho(f)(e^i) - (i_{\varphi(f)}\otimes\id)(\nabla(e^i))\big] \Big|& \textrm{since }\varepsilon^*_i\circ j_M = F\circ j_A\circ e^*_i,\\
	&= \sum_i \Big| \tilde\Delta\circ e^*_i\big[\rho(f)(e^i) - (i_{\varphi(f)}\otimes\id)(\nabla(e^i))\big] \Big|& \textrm{since }F\circ j_A = \tilde\Delta,\\
	&= \tilde\Delta(\Div^{(\nabla,\varphi,\rho)}(f))\,.
\end{align*}
This concludes the proof.\qed

\begin{remark}
We can also go in the opposite direction. Namely, we have the functor
\[
	\Xi_A = \,\cdot\,\otimes_A\mathbb{K}: A^\mathrm{e}\textrm{-}\textbf{Mod}\to A\textrm{-}\textbf{Mod}
\]
and the natural transformation $j'\colon \id_{A^\mathrm{e}\textrm{-}\textbf{Mod}} \to \Xi_A$ given by
\[
	j'_N\colon N\to N\otimes_A\mathbb{K}:n\mapsto n\otimes 1\,,
\]
which is a module homomorphism over $\id\otimes \,\varepsilon$. We have the natural isomorphism $\Xi_A\Phi_A\cong \id$, but not the other way around.\\
\end{remark}

Now, we move on to Lie algebras. Let $\mathfrak{g}$ be a Lie algebra over $\mathbb{K}$ and $U\mathfrak{g}$ its enveloping algebra, which is endowed with the Hopf algebra structure given, for $x\in\mathfrak g$, by
\[
	\Delta(x) = x\otimes 1 + 1\otimes x,\;\; S(x) = -x, \;\textrm{ and } \; \varepsilon(x) = 0.
\]

\begin{definition}\ 
\begin{itemize}
	\item For a $\mathfrak g$-module $M$, $\Der_\mathrm{Lie}(\mathfrak g,M)$ is the space of all $\mathbb{K}$-linear derivations (i.e., Lie algebra $1$-cocycles) on $\mathfrak g$ into $M$. Set $\Der_\mathrm{Lie}(\mathfrak g) = \Der_\mathrm{Lie}(\mathfrak g,\mathfrak g)$.
	\item For $f\in\Der_\mathrm{Lie}(\mathfrak g)$, denote by $U\!f$ the natural extension of $f$ to $U\mathfrak g$.
	\item $\Omega^1_\mathrm{Lie}\mathfrak{g} = \Ker(\varepsilon\colon U\mathfrak g \to \mathbb{K})$ is the space of \textit{Lie $1$-forms}. This is a left $U\mathfrak g$-submodule of $U\mathfrak g$, and we have the universal derivation $d\colon \mathfrak{g}\to\Omega^1_\mathrm{Lie}\mathfrak{g}$ given by $dx = x$. The pair $(\Omega^1_\mathrm{Lie}\mathfrak{g},d)$ is characterised by the following universal property: for any $\mathfrak{g}$-module $M$ and a derivation $f\colon\mathfrak{g}\to M$, there is a unique $U\mathfrak{g}$-module homomorphism $i_f\colon \Omega^1_\mathrm{Lie}\mathfrak{g} \to M$ such that $i_f\circ d = f$ holds.
	\item Let $a\in U\mathfrak g$ and $x,y\in \mathfrak g$. The element $ax$ of $\Omega^1_\mathrm{Lie}\mathfrak{g}$ is formally written as $adx$, with the Leibniz rule
	 \[
	 	d[x,y] = xdy - ydx.
	\]
	\item  For $f\in\Der_\mathrm{Lie}(\mathfrak g)$, the \textit{Lie derivative} $L_f \colon \Omega^1_\mathrm{Lie}\mathfrak{g}\to \Omega^1_\mathrm{Lie}\mathfrak{g}$ is defined by the formula
	\[
		L_f(adx) = (U\!f)(a)dx + adf(x) \quad \textrm{for } a\in U\mathfrak g\textrm{ and } x\in\mathfrak g,.
	\]
	\item A Lie algebra $\mathfrak g$ is said to be \textit{formally smooth} if $\mathfrak g$ is finitely generated as a Lie algebra, and $\Omega^1_\mathrm{Lie}\mathfrak{g}$ is dualisable as a $U\mathfrak g$-module.
\end{itemize}
\end{definition}

\begin{lemma}\label{lem:lie1form}
We have the isomorphism of $U\mathfrak g$-bimodules $\Phi_{U\mathfrak g}(\Omega^1_\mathrm{Lie}\mathfrak g) \cong \Omega^1U\mathfrak g$ given by 
\[
	adx\otimes b \mapsto -\bar b\tilde\Delta(a)\cdot dx\;\textrm{ for } a\in U\mathfrak g \textrm{ and } x\in \mathfrak g\,.
\]
\end{lemma}
\noindent Proof. By the definition of Lie $1$-forms, we have an exact sequence
\[\begin{tikzcd}[cramped, ampersand replacement = \&]
	0\arrow[r] \&  \Omega^1_\mathrm{Lie}\mathfrak g \arrow[r] \& U\mathfrak g \arrow[r,"\varepsilon"]\& \mathbb{K} \arrow[r] \& 0
\end{tikzcd}\]
of left $U\mathfrak g$-modules. Applying the functor $\Phi_{U\mathfrak g}$, the isomorphisms (\ref{eq:bimodisom}) yields a commutative diagram
\[\begin{tikzcd}[cramped, ampersand replacement = \&]
	0\arrow[r] \&  \Phi_{U\mathfrak g}(\Omega^1_\mathrm{Lie}\mathfrak g) \arrow[r] \& \Phi_{U\mathfrak g}(U\mathfrak g) \arrow[r,"\Phi_A(\varepsilon)"]\arrow[d,"\sim"]\& \Phi_{U\mathfrak g}(\mathbb{K}) \arrow[r]\arrow[d,"\sim"] \& 0\\
	0\arrow[r]\&\Omega^1 U\mathfrak g\arrow[r]\& U\mathfrak g^{\otimes 2} \arrow[r,"\mu"]\& U\mathfrak g \arrow[r]\& 0
\end{tikzcd}\]
of $U\mathfrak g^\mathrm{e}$-modules with the two rows exact. Hence, we obtain the isomorphism $\Phi_{U\mathfrak g}(\Omega^1_\mathrm{Lie}\mathfrak g) \cong \Omega^1U\mathfrak g$. The map is explicitly given by $adx\otimes b \mapsto -\bar b\tilde\Delta(a)\cdot dx$ for $a,b\in U\mathfrak g$ and $x\in \mathfrak g$\,, by chasing the following diagram:
 \[\begin{tikzcd}[cramped, ampersand replacement = \&]
 	adx\otimes b \arrow[r,mapsto] \& ax\otimes b\arrow[d,mapsto]\\
	-\bar b\tilde\Delta(a)\cdot dx\arrow[r,mapsto]\& \bar b\tilde\Delta(a)\cdot (x\otimes 1 - 1\otimes x) = \bar b\tilde\Delta(ax)\,.
\end{tikzcd}\]
This completes the proof.\qed

\begin{remark}\label{rem:liesmooth}
Formally smooth Lie algebras are of cohomological dimension one (see Proposition 21.1.19 of \cite{Ginzburg:aa}), and one example is a finitely generated free Lie algebra $L(z_1,\dotsc, z_r)$. The existence of cohomologically one-dimensional non-free Lie algebras is an open problem over a field of characteristic zero; see the Introduction of \cite{ZUSMANOVICH2019288} and references there. If $\mathfrak g$ is formally smooth, $U\mathfrak g$ is automatically formally smooth as an associative algebra as the functor (\ref{eq:hopfadj}) admits the right adjoint $\Psi_A$, as seen in Remark \ref{rem:adj}.
\end{remark}

Now suppose that $\mathfrak g$ is formally smooth and that we are given a connection $\nabla\colon \Omega^1_\mathrm{Lie}\mathfrak{g}\to \Omega^1U\mathfrak{g}\otimes_{U\mathfrak{g}}\Omega^1_\mathrm{Lie}\mathfrak{g}$\,. We consider the divergence map, defined in Definition \ref{def:conn}, associated with the following data: $B = U\mathfrak g$, $M = \Omega^1_\mathrm{Lie}\mathfrak{g}$, $\mathfrak d = \Der_\mathrm{Lie}(\mathfrak g)$,
\begin{align*}
	\varphi_\mathrm{Lie}&\colon \Der_\mathrm{Lie}(\mathfrak g) \to \Der_\mathbb{K}(U\mathfrak g)\hspace{12pt}\colon f \mapsto  U\!f,\textrm{ and}\\
	\rho_\mathrm{Lie}&\colon \Der_\mathrm{Lie}(\mathfrak g) \to \End_\mathbb{K}(\Omega^1_\mathrm{Lie}\mathfrak g)\colon f \mapsto L_f.
\end{align*}
From these, we obtain the associated divergence map
\[
	\sdiv^{\nabla}\colon \Der_\mathrm{Lie}(\mathfrak g) \to |U\mathfrak g|.
\]
This \textit{single} divergence map is related to the double ones in the associative setting:

\begin{corollary}
Let $U\nabla$ be the connection on $\Omega^1U\mathfrak g$ obtained from $\Phi_{U\mathfrak g}(\nabla)$ via the isomorphism in Lemma \ref{lem:lie1form}. Then, we have
\[
	\Div^{U\nabla}(U\!f) = \tilde\Delta(\sdiv^{\nabla}(f))
\]
for $f\in\Der_\mathrm{Lie}(\mathfrak g)$ in $U\mathfrak g^\mathrm{e}$.
\end{corollary}
\noindent Proof. Recall that the double divergence is associated with the derivation action $(\Der_\mathbb{K}(U\mathfrak g),\varphi_\mathrm{assoc},\rho_\mathrm{assoc})$. Therefore, $\Div^{\nabla}\circ \,\varphi_\mathrm{Lie}$ is associated with $(\Der_\mathrm{Lie}(\mathfrak g), \varphi_\mathrm{assoc}\circ\,\varphi_\mathrm{Lie}, \rho_\mathrm{assoc}\circ\,\varphi_\mathrm{Lie})$. \vspace{2pt}On the other hand, we have $\Imag\varphi = \Der_\mathrm{Hopf}(U\mathfrak g)$. Then, by Theorem \ref{thm:hopfconn}, $\tilde\Delta\circ\sdiv^{\nabla}$ is equal to $\Div^{(\Phi_{U\mathfrak g}(\nabla),\Phi_{U\mathfrak g}(\varphi_\mathrm{Lie}),\Phi_{U\mathfrak g}(\rho_\mathrm{Lie}))}$. We will show that these two are equal.

By definition, $\Phi_{U\mathfrak g}(\varphi_\mathrm{Lie}) = \varphi_\mathrm{assoc}\circ\,\varphi_\mathrm{Lie}$. Next, we check that the following diagram is commutative for any $f\in\Der_\mathrm{Lie}(\mathfrak g)$:
\[\begin{tikzcd}[ampersand replacement = \&]
	\Phi_{U\mathfrak g}(\Omega^1_\mathrm{Lie}\mathfrak g) \arrow[r,"\Phi_{U\mathfrak g}(\rho_\mathrm{Lie})(f)"]\arrow[d,"\sim"] \&\Phi_{U\mathfrak g}(\Omega^1_\mathrm{Lie}\mathfrak g)\arrow[d,"\sim"] \& adx\otimes b \arrow[r,mapsto] \arrow[d,mapsto]\& {U\!f}(a)dx\otimes b + adf(x)\otimes b + adx\otimes U\!f(b)\arrow[d,mapsto] \\
	\Omega^1U\mathfrak g \arrow[r,"\rho_\mathrm{assoc}\circ\,\varphi_\mathrm{Lie}(f)"]\& \Omega^1U\mathfrak g \& -\bar b\tilde\Delta(a)\cdot dx \arrow[r,mapsto]\& -\bar b\tilde\Delta(U\!f(a))\cdot dx -\bar b\tilde\Delta(a)\cdot df(x) - \overline{U\!f(b)}\tilde\Delta(a)\cdot dx
\end{tikzcd}\]
This shows that the two derivation actions agree via the isomorphism in \ref{lem:lie1form}. Hence we have
\[
	\tilde\Delta\circ\sdiv^{\nabla} = \Div^{(\Phi_{U\mathfrak g}(\nabla),\Phi_{U\mathfrak g}(\varphi_\mathrm{Lie}),\Phi_{U\mathfrak g}(\rho_\mathrm{Lie}))} = \Div^{(U\nabla,\varphi_\mathrm{assoc}\circ\,\varphi_\mathrm{Lie}, \rho_\mathrm{assoc}\circ\,\varphi_\mathrm{Lie})} = \Div^{U\nabla}\circ \,\varphi_\mathrm{Lie}\,,
\]
which completes the proof. \qed

\begin{remark}
All the constructions above work just as well on topological algebras, provided everything is continuous.
\end{remark}

Let us get back to the flat connections in the previous section.

\begin{example}\ 
\begin{enumerate}[(1)]
	\item The connection $\nabla\!_z$ on $\Omega^1\mathbb{K}\langle z_1,\dotsc,z_r\rangle$ can be written as $U\nabla'\!\!_z$ where the connection $\nabla'\!\!_z$ on $\Omega^1_\mathrm{Lie}L(z_1,\dotsc, z_r)$ is defined by $\nabla'\!\!_z(dz_i)=0$ for all $i$.
	\item Similarly, the connection $\nabla\!_\mathcal{C}$ on $\Omega^1\mathbb{K}F_r$ is induced from $\Phi_{\mathbb{K}F_r}(\nabla'\!\!_\mathcal{C})$, where the connection $\nabla'\!\!_\mathcal{C}$ on
	\[
		\Ker(\varepsilon) = \bigoplus_{c\in \mathcal{C}}A\otimes\mathbb{K}\cdot(c-1)
	\]
	is defined by $\nabla'\!\!_\mathcal{C}(1\otimes (c-1))=0$ for all $c\in\mathcal{C}$. We use the isomorphism $\Phi_{\mathbb{K}F_r}(\Ker(\varepsilon))\cong \Omega^1\mathbb{K}F_r$ analogous to Lemma \ref{lem:lie1form} to transport $\Phi_{\mathbb{K}F_r}(\nabla'\!\!_\mathcal{C})$ onto $\Omega^1\mathbb{K}F_r$.
\end{enumerate}
\end{example}
\noindent On the other hand, $\nabla\!_\mathcal{C}$ is not induced from a connection on a Lie algebra since $\mathbb{K}F_r$ is not of the form $U\mathfrak g$ for any $\mathfrak g$. However, we have the following:

\begin{proposition}
Let $\widehat{\mathbb{K}F_{2g+n}}$ be the completion of $\mathbb{K}F_{2g+n}$ with respect to the augmentation ideal
\begin{align*}
	\Ker(\varepsilon\colon \mathbb{K}F_{2g+n}\to\mathbb{K})\,.
\end{align*}
Then, the continuous extension of the connection $\nabla\!_\mathcal{C}$ is induced from a connection on $\Omega^1_\mathrm{Lie}\mathbb{P}(\widehat{\mathbb{K}F_{2g+n}})$, where $\mathbb{P}$ denotes the primitive part.
\end{proposition}
\noindent Proof. After the completion, $\widehat{\mathbb{K}F_{2g+n}}$ is isomorphic, as a complete Hopf algebra, to the completed free associative algebra $\mathbb{K}\langle\!\langle z_1,\dotsc,z_{2g+n}\rangle\!\rangle$ by the following map:
\[
	\widehat{\mathbb{K}F_{2g+n}}\cong \mathbb{K}\langle\!\langle z_1,\dotsc,z_{2g+n}\rangle\!\rangle\colon x_i \mapsto e^{z_i}\;\textrm{ for } x_i\in \mathcal{C}\,,
\]
which induces an isomorphism on the primitive parts
\[	
	\mathbb{P}(\widehat{\mathbb{K}F_{2g+n}})\cong L(\!(z_1,\dotsc,z_{2g+n})\!)
\]
where $L(\!(z_1,\dotsc,z_{2g+n})\!)$ is the completed free Lie algebra of rank $2g+n$. Under this identification, we have
\begin{align*}
	dx_i\,x_i^{-1} &= 1\otimes 1 - x_i \otimes x_i^{-1} = \tilde\Delta(1 - x_i) =  \tilde\Delta(1 - e^{z_i})\\
	&= \tilde\Delta\left(\frac{e^{z_i}-1}{z_i}\right)\tilde\Delta(-z_i) = \tilde\Delta\left(\frac{e^{z_i}-1}{z_i}\right)\cdot dz_i = \tilde\Delta\left(\frac{e^{z_i}-1}{z_i}\,dz_i\right)\,
\end{align*}
which shows that the continuous extension of $\nabla\!_\mathcal{C}$ is induced from the connection $\nabla''_\mathcal{C}$ on $\Omega^1_\mathrm{Lie}L(\!(z_1,\dotsc,z_{2g+n})\!)$ defined by
\[
	\nabla''\!\!_\mathcal{C}\left(\frac{e^{z_i}-1}{z_i}\,dz_i\right)=0.
\]
This completes the proof.\qed\\

Finally, recall the Satoh trace map introduced in \cite{SATOH2012709}. Let $H = \mathbb{K}\{z_1,\dotsc,z_n\}$ be a free $\mathbb{K}$-module, $H^* = \Hom_\mathbb{K}(H,\mathbb{K})$ its dual space, and $\mathfrak g = L(z_1,\dotsc,z_n)$. Then $A = U\mathfrak g$ is isomorphic to $\mathbb{K}\langle z_1,\dotsc,z_n\rangle$ as a Hopf algebra. The Satoh trace map $\Tr_\mathrm{Satoh}$ is defined as the composite
\[\begin{tikzcd}[ampersand replacement = \&]
	\Tr_\mathrm{Satoh}\colon H^*\otimes L(H) \arrow[r, "\mathrm{incl}"] \&  H^*\otimes A \arrow[r,"\mathrm{cont}"] \& A \arrow[r, "\mathrm{proj}"]\& \text{\textbar}A\text{\textbar},
\end{tikzcd}\]
where the contraction map is defined by
\[
	\mathrm{cont}(z^*_j\otimes z_{i_1}\cdots z_{i_\ell}) =  \delta_{i_1,j}\cdot z_{i_2}\cdots z_{i_\ell}\,.
\]
Satoh showed in \cite{SATOH2012709} that the kernel of $\Tr_\mathrm{Satoh}$ stably coincides with the image of the Johnson homomorphism on the automorphism group of a free group. On the other hand, the map $\Tr_\mathrm{Satoh}$ is exactly the map $\sdiv^{\nabla'\!\!_z}$, as is also observed in Section 9 of \cite{akkn}. Our construction enables us to deduce the important $1$-cocycle property of the Satoh trace with no combinatorial difficulties. \\

For an introduction into geometry over an operad, see the last section of \cite{Ginzburg:aa}, for example. The divergence maps on algebras over other operads can be constructed similarly as long as universal enveloping algebras are defined. Some aspects in this direction are also discussed in Appendix B in \cite{Powell:2021aa}.\\

\section{Perfect Complexes and Homological Connections}\label{sec:perf}

In this section, we define the divergence associated with a connection on a not-necessarily-projective module. Since the existence of a connection is closely tied in with projectivity, as we have seen in Remark \ref{rem:connproj}, we have to modify the definition to work with them. \\

\noindent\textbf{Conventions.} A projective resolution of a $B$-module $M$ is of the form $P = (\,\cdots\to P_1\xrightarrow{\partial_1}P_0 \xrightarrow{\partial_0}M\to 0\,)$. We set $P_{-1}=M$.

\begin{definition}
For a $B$-module $M$ and $f\in\Der_\mathbb{K}(B)$, a $\mathbb{K}$-linear map $u:M\to M$ is called an \textit{$f$-derivation} if
\[
	u(bm) = f(b)m + bu(m)
\]
holds for $b\in B$ and $n\in M$. Denote by $\Der_f(M)$ the space of $f$-derivations on $M$.
\end{definition}

\begin{lemma}\label{lem:fder}
Let $f,g\in\Der_\mathbb{K}(B)$.
\begin{enumerate}[(1)]
	\item If $\Der_f(M)$ is non-empty, it is an affine space modelled on $\End_B(M)$.
	\item If $M$ is projective, $\Der_f(M)$ is non-empty.
	\item If $u\in\Der_f(M)$ and $u'\in\Der_g(M)$, then $[u,u']$ is an $[f,g]$-derivation. \label{lem:item:comm}
\end{enumerate}
\end{lemma}
\noindent Proof. (1) If $u, u'\in\Der_f(M)$, then $u-u'\colon M\to M$ is clearly $B$-linear.\\
(2) Assume that $M$ is projective so that we can take a surjective $B$-module map $q\colon B^{\oplus I} \twoheadrightarrow M$ for some index set $I$, as well as its $B$-linear section $s\colon M\to B^{\oplus I}$. Define $\gamma:M\to M$ by the composite
\[
	\gamma\colon M\xrightarrow{s} B^{\oplus I} \xrightarrow{f^{\oplus I} } B^{\oplus I} \xrightarrow{q} M\,.
\]
Then, we have, for $b\in B$ and $m\in M$,
\begin{align*}
	\gamma(bm) &= q\circ f^{\oplus I} \circ s(bm)\\
	&= q\circ f^{\oplus I}(b s(m))\\
	&= q(f(b) s(m) + bf^{\oplus I}(s(m)))\\
	& = f(b)m + b\gamma(m)\,,
\end{align*}
using the fact that $f$ is a derivation and $q$ and $s$ are $B$-linear. This shows that $\gamma$ is an $f$-derivation.\\
(3) Let $b\in B$ and $m\in M$. We proceed by direct computation:
\begin{align*}
	[u,u'](bm) &= u(g(b)m + bu'(m)) - u'(f(b)m + bu(m))\\
	&= f(g(b))m + g(b)u(m) + f(b)u'(m) + bu(u'(m))\\
	&\hspace{40pt} - g(f(b))m - f(b)u'(m) - g(b)u(m) - bu'(u(m))\\
	&= f(g(b))m + bu(u'(m)) - g(f(b))m  - bu'(u(m))\\
	&= [f,g](b)m + b[u,u'](m)\,.
\end{align*}
This completes the proof.\qed

\begin{remark}
Given a derivation action $(\mathfrak d,\varphi,\rho)$ on $M$ considered in Definition \ref{def:deract}, the compatibility condition says that, for $f\in\mathfrak d$, each $\rho(f)$ is a $\varphi(f)$-derivation.\\
\end{remark}

\begin{proposition}\label{prop:derlift}
Let $(P,\partial)$ be a projective resolution of $M$, $f\in\Der_\mathbb{K}(B)$, and $u\colon M\to M$ an $f$-derivation. Then, $u$ admits a lift $\,\lambda[u]:P\to P$ comprised of $f$-derivations satisfying $[\partial,\lambda[u]]=0$. Furthermore, such a lift is unique up to homotopy.
\end{proposition}
\noindent Proof. Take an arbitrary family of $f$-dervations $(\gamma_n:P_n\to P_n)_{n\geq 0}$ by Lemma \ref{lem:fder}. We construct maps by following steps: first, we set
\begin{align*}
	\lambda[u]_{-1} &=u :M \to M,
	\intertext{and, for $n\geq0$,}
	\varphi[u]_n &= \lambda[u]_{n-1}\circ\partial_n - \partial_n\circ \gamma_n:P_n\to P_{n-1},\\
	\tilde\varphi[u]_n&: P_n\to P_n\:\textrm{ so that }\: \partial_n\circ\tilde\varphi[u]_n = \varphi[u]_n ,\textrm{ and}\\
	\lambda[u]_n &= \tilde\varphi[u]_n + \gamma_n:P_n\to P_n.
\end{align*}
We check that this inductive procedure is well-defined. Let $n\geq 0$ and suppose that $\lambda[u]_{n-1}$ is an $f$-derivation. Then, $\varphi[u]_n$ is a $B$-module map: for $b\in B$ and $p\in P_n$,
\begin{align*}
	\varphi[u]_n(bp) &= \lambda[u]_{n-1}(\partial_n(bp)) - \partial_n(\gamma_n(bp))\\
	&= \lambda[u]_{n-1}(b\partial_n(p)) - \partial_n(f(b)p + b\gamma_n(p))\\
	&= f(b)\partial_n(p) + b\lambda[u]_{n-1}(\partial_n(p)) - f(b)\partial_n(p) + b\partial_n(\gamma_n(p))\\
	&= b\varphi[u]_n(p).
\end{align*}
Next, we check $\partial_{n-1}\circ\varphi[u]_n=0$. If $n=0$, this is true since $\partial_{-1}$ is defined to be zero. If $n\geq1$, we have
\begin{align*}
	\partial_{n-1}\circ\varphi[u]_n &= \partial_{n-1}\circ(\lambda[u]_{n-1}\circ\partial_n - \partial_n\circ \gamma_n)\\
	&= \partial_{n-1}\circ(\tilde\varphi[u]_{n-1} + \gamma_{n-1})\circ\partial_n\\
	&= (\varphi[u]_{n-1} +  \partial_{n-1}\circ \gamma_{n-1})\circ\partial_n\\
	&= (\lambda[u]_{n-2}\circ\partial_{n-1} - \partial_{n-1}\circ \gamma_{n-1} +  \partial_{n-1}\circ \gamma_{n-1})\circ\partial_n\\
	&= 0,
\end{align*}
which shows that we can take a lift of $\varphi[u]_n$. Lastly, it is clear that $\lambda[u]_n$ is an $f$-derivation. This shows the well-definedness. Next, we have
\begin{align*}
	\partial_n\circ\lambda[u]_n &= \partial_n\circ(\tilde\varphi[u]_n + \gamma_n)\\
	&= \varphi[u]_n + \partial_n\circ \gamma_n\\
	&= \lambda[u]_{n-1}\circ\partial_n - \partial_n\circ \gamma_n + \partial_n\circ \gamma_n\\
	&= \lambda[u]_{n-1}\circ\partial_n,
\end{align*}
which is briefly denoted by $[\partial,\lambda[u]]=0$.

Now suppose we have two such lifts $\lambda[u]$ and $\lambda[u]'$. Then the difference $\lambda[u]-\lambda[u]'$ is a $B$-module map by Lemma \ref{lem:fder}, and at the same time, is a lift of the zero map $0:M\to M$. Hence, this is null-homotopic, and this completes the proof.\qed\\

\begin{definition}
Let $(C,\partial)$ be a chain complex of $B$-modules.
\begin{itemize}
	\item We put $\IHom^0_B(C,C) = \{(\psi_n:C_n\to C_n)_n:B\textrm{-module maps}\}$, the space of degree zero maps. $\IHom$ stands for the internal hom-set.
	\item $(C,\partial)$ is called a \textit{perfect complex} if it is of finite length and each $C_n$ is $B$-dualisable. In this case, the trace map is defined by
	\begin{align*}
		\Tr:\IHom^0_B(C,C)\to|B|: (\psi_n)_n\mapsto \sum_n(-1)^n\Tr(\psi_n).
	\end{align*}
	\item A $B$-module $M$ is \textit{perfect} if it admits a projective resolution by a perfect complex.
\end{itemize}
\end{definition}

\begin{lemma}\label{lem:tracehomotopy}
Let $(C,\partial)$ be a perfect complex of $B$-modules and $(h_n:C_n\to C_{n+1})_n$ be a family of $B$-module maps. Then we have $\Tr([\partial,h])=0$.
\end{lemma}
\noindent Proof. By definition, we have
\begin{align*}
	\Tr([\partial,h]) &= \sum_n(-1)^n \Tr(\partial_{n+1}\circ h_n + h_{n-1}\circ\partial_n)\\
	&= \sum_n (-1)^n \Tr(\partial_{n+1}\circ h_n) + \sum_n(-1)^{n+1}\Tr(h_n\circ\partial_{n+1}).
\end{align*}
Since $\Tr(h_n\circ\partial_{n+1}) = \Tr(\partial_{n+1}\circ h_n)$ by Proposition \ref{prop:HTtraceLie}, we obtain $\Tr([\partial,h])=0$.\qed\\

Now, we can define a connection on a $B$-module and the divergence associated with it. 
\begin{definition}\label{def:homconn}
Let $M$ be a $B$-module. 
\begin{itemize}
	\item A \textit{homological connection} on $M$ is a pair $\nabla = (P, \nabla_\bullet)$ where $P$ is a projective resolution of $M$ and\\
	$\{\nabla_n:P_n\to \Omega^1B\otimes_BP_n\}_{n\geq0}$ is a family of usual connections.
	\item The \textit{curvature} of $\nabla$ is defined to be $R = \{(\nabla_n)^2:P_n\to \Omega^2B\otimes_BP_n\}_{n\geq0}$, which is a family of $B$-module maps. A homological connection $\nabla$ is \textit{flat} if the curvature $R$ is the zero map.
\end{itemize}
\end{definition}
\noindent Note that the collection $\nabla$ of connections above is \textit{not} required to be a chain map.

\begin{definition}
Let $\nabla = (P,\nabla_\bullet)$ be a homological connection on $M$ with $P$ perfect, and $(\mathfrak d,\varphi,\rho)$ be a derivation action on $M$. The associated \textit{divergence} $\Div^{(\nabla,\varphi,\rho)}$ is defined by
\[
	\Div^{(\nabla,\varphi,\rho)}\!:\mathfrak d\to |B|: f\mapsto  \sum_{n\geq0}(-1)^n\Tr\left(\lambda[\rho(f)]_n - (i_{\varphi(f)}\otimes \id)\circ\nabla_n\right).
\]
Note that $\lambda[\rho(f)]_n$ and $(i_{\varphi(f)}\otimes \id)\circ\nabla_n$ are both $\varphi(f)$-derivations, so that we can take the trace of their difference. This is well-defined by Proposition \ref{prop:derlift} and Lemma \ref{lem:tracehomotopy}.
\end{definition}

\begin{proposition}
If $\nabla = (P,\nabla_\bullet)$ is flat, $\Div^{(\nabla,\varphi,\rho)}$ is a Lie algebra $1$-cocycle.
\end{proposition}
\noindent Proof. Let $f,g\in\mathfrak d$ and take lifts $\lambda[\rho(f)]$ and $\lambda[\rho(g)]$ of $\rho(f)$ and $\rho(g)$ to $P$, respectively. Then $[\lambda[\rho(f)],\lambda[\rho(g)]]$ is a chain map, and also a $[f,g]$-derivation by Lemma \ref{lem:fder} \ref{lem:item:comm}. Therefore $\lambda[\rho([f,g])]$ is homotopic to $[\lambda[\rho(f)],\lambda[\rho(g)]]$ by Proposition \ref{prop:derlift}. The rest is analogous to the proof of Proposition \ref{prop:flatcocyc}.\qed\\

We also have a relaxed version of smoothness (see Definition 8.1.2 of \cite{Kontsevich_2008}).

\begin{definition}
An algebra $A$ is said to be \textit{homologically smooth} if $A$ is finitely generated as a $\mathbb{K}$-algebra and $\Omega^1A$ is a perfect module over $A^\mathrm{e}$.
\end{definition}

In this case, we obtain a divergence map $\Div^\nabla\colon\Der_\mathbb{K}(A)\to|A^\mathrm{e}|$, associated with $\nabla$ and the default action given in (\ref{eq:defaultaction}), generalising the construction given in Section \ref{sec:condiv}.\\

\section{The Closed Surface Case}\label{sec:closed}
Let $\Sigma_{g,1}$ be an oriented surface of genus $g$ with one boundary component, and pick a base point on the boundary. Denote by $\zeta = \zeta_0$ the only simple boundary loop in $\pi_1(\Sigma_{g,1})$. Recall that we fixed an isomorphism
\[
	\pi_1(\Sigma_{g,1})\cong F_{2g} = \langle a_i,b_i\,(1\leq i\leq g)\rangle
\]
in Section \ref{sec:cob} so that $\zeta = (a_1,b_1)\cdots(a_g,b_g)$ holds. Then the closed surface $\Sigma_g$ is obtained by capping the boundary with a disk, with the base point induced from $\Sigma_{g,1}$. Put $\pi = \pi_1(\Sigma_g)$.

\begin{definition}
The map $v\colon|\mathbb{K}\pi|\to \mathrm{HH}^1(\mathbb{K}\pi)$ is defined by
\[
	v(\alpha)(x) = \sum_{p\in\alpha\cap x} \mathrm{sign}(\alpha,x;p) \,\alpha\ast_px
\]
for generic representatives of a free loop $\alpha$ and $x\in\pi$. It is originally denoted by $\rho$ in \cite{Vaintrob2007}, and its well-definedness is due to Vaintrob (Lemma 1 of \cite{Vaintrob2007}).
\end{definition}

In this section, we will show the following theorem, which gives an algebraic description of the Turaev cobracket on a closed surface. Put $R = \mathbb{K}\langle\zeta\rangle$, the subalgebra of $\mathbb{K}F_{2g}$ generated by $\zeta$. Then, we have the natural map $\Der_R(\mathbb{K}F_{2g})\to \Der_\mathbb{K}(\mathbb{K}\pi)$ since $\pi$ is the quotient of $F_{2g}$ by the normal subgroup generated by $\zeta$.

\begin{theorem}\label{thm:cob}
Let $\mathcal{C} = (a_i,b_i)_{1\leq i\leq g}$ be a free-generating system of $\pi_1(\Sigma_{g,1})$ with $\zeta = (a_1,b_1)\cdots(a_g,b_g)$. Then we have the commutative diagram
\begin{align}\label{eq:twosquare}
	\begin{tikzcd}[ampersand replacement = \&]
	\text{\textbar}\mathbb{K}F_{2g}\text{\textbar} \arrow[r, "\sigma"] \arrow[dd,two heads] \& \Der_R(\mathbb{K}F_{2g}) \arrow[r,"\Div^\mathcal{C}"]\arrow[d] \& \text{\textbar}\mathbb{K}F_{2g}\text{\textbar}^{\otimes 2}\arrow[d,two heads]   \phantom{.}\\
	\& \Der_\mathbb{K}(\mathbb{K}\pi) \arrow[d,two heads] \& \text{\textbar}\mathbb{K}\pi\text{\textbar}^{\otimes 2} \arrow[d,two heads] \\
	\text{\textbar}\mathbb{K}\pi\text{\textbar} \arrow[r, "v"] \& \mathrm{HH}^1(\mathbb{K}\pi) \arrow[r,"\Div^{\nabla'}"] \& \text{\textbar}\mathbb{K}\pi/\mathbb{K}1\text{\textbar}^{\otimes 2} .
	\end{tikzcd}
\end{align}
for some homological connection $\nabla'$ on $\Omega^1\mathbb{K}\pi$.
\end{theorem}

\begin{corollary}\label{cor:cob}
The composite $\Div^{\nabla'}\!\!\circ \,v$ is equal to the Turaev cobracket $\delta$.
\end{corollary}
\noindent Proof. By Theorem \ref{thm:akkn}, $\Div^\mathcal{C}\circ\,\sigma$ is equal to the suitably framed version of the Turaev cobracket $\delta^\mathsf{fr}$. On the other hand, we have the following commutative diagram
\[\begin{tikzcd}[ampersand replacement = \&]
	\text{\textbar}\mathbb{K}F_{2g}\text{\textbar} \arrow[r, "\delta^\mathsf{fr}"]\arrow[d, two heads]  \& \text{\textbar}\mathbb{K}F_{2g}\text{\textbar}^{\otimes 2}\arrow[d,two heads]\\
	\text{\textbar}\mathbb{K}\pi\text{\textbar} \arrow[r, "\delta"]  \& \text{\textbar}\mathbb{K}\pi/\mathbb{K}1\text{\textbar}^{\otimes 2}
\end{tikzcd}\]
with the same maps on the top, left, and right with the outer square of (\ref{eq:twosquare}). Therefore, the composite $\Div^{\nabla'}\!\!\circ\, v$ coincides with the Turaev cobracket $\delta$ by the surjectivity of the map $|\mathbb{K}F_{2g}|\to |\mathbb{K}\pi|$.\qed\\

To prove Theorem \ref{thm:cob}, we first construct a homological connection $\nabla'$. To do so, we have to choose a projective resolution of $\mathbb{K}\pi$, which we can take to be a perfect complex since $\mathbb{K}\pi$ is homologically smooth.

\begin{lemma}
Let $V = \mathbb{K}\{a_i,b_i\}_{1\leq i\leq g}$ be a free $\mathbb{K}$-module of rank $2g$. We have a $\mathbb{K}\pi^\mathrm{e}$-free resolution of $\mathbb{K}\pi$:
\begin{align}
	\begin{tikzcd}[cramped,ampersand replacement = \&]
	0 \arrow[r] \&\mathbb{K}\pi\otimes\mathbb{K}\pi \arrow[r, "d_1"] \& \mathbb{K}\pi\otimes V\otimes \mathbb{K}\pi \arrow[r,"d_0"] \& \mathbb{K}\pi\otimes\mathbb{K}\pi \arrow[r,"\mu"] \& \mathbb{K}\pi\arrow[r] \& 0\,,
	\end{tikzcd}\label{eq:resclean}
\end{align}\vspace{-20pt}
\begin{align*}
	d_1(1\otimes 1) &= \sum_{1\leq i\leq g} (a_1,b_1)\cdots(a_{i-1},b_{i-1}) \Big( 1\otimes a_i\otimes b_ia_i^{-1}b_i^{-1} + a_i\otimes b_i\otimes a_i^{-1}b_i^{-1}\\
	&\hspace{100pt}- a_ib_ia_i^{-1}\otimes a_i\otimes a_i^{-1}b_i^{-1} - a_ib_ia_i^{-1}b_i^{-1}\otimes b_i\otimes b_i^{-1} \Big)(a_{i+1},b_{i+1})\cdots (a_g,b_g),\\
	d_0(1\otimes v \otimes 1) &= 1\otimes v - v\otimes 1.
\end{align*}
\end{lemma}
\noindent Proof. We follow the method by R. Lyndon \cite{Lyndon50}. By Section 11 of \cite{Lyndon50}, we have a left $\mathbb{K}\pi$-free resolution of $\mathbb{K}$ using the Fox derivative:
\begin{align*}
	\begin{tikzcd}[cramped,ampersand replacement = \&]
	0 \arrow[r] \& \mathbb{K}\pi \arrow[r, "\partial_1"] \& \mathbb{K}\pi\otimes V \arrow[r,"\partial_0"] \& \mathbb{K}\pi \arrow[r,"\varepsilon"] \& \mathbb{K}\arrow[r] \& 0\,,
	\end{tikzcd}
\end{align*}\vspace{-20pt}
\begin{align*}
	\partial_1(1) &=  \sum_{1\leq i\leq g} \left( \frac{\partial \zeta}{\partial a_i}\otimes a_i + \frac{\partial \zeta}{\partial b_i}\otimes b_i \right)\\
	&= \sum_{1\leq i\leq g} (a_1,b_1)\cdots(a_{i-1},b_{i-1})\left( 1\otimes a_i + a_i\otimes b_i - a_ib_ia_i^{-1}\otimes a_i - a_ib_ia_i^{-1}b_i^{-1}\otimes b_i \right),\\
	\partial_0(1\otimes c) &= c-1 \quad \textrm{for } c\in\{a_i,b_i\}_{1\leq i\leq g}\,.
\end{align*}
Applying the functor (\ref{eq:hopfadj}) and isomorphisms in (\ref{eq:bimodisom}), we obtain a left $\mathbb{K}\pi^\mathrm{e}$-free resolution
\begin{align*}
	\begin{tikzcd}[cramped,ampersand replacement = \&]
	0 \arrow[r] \& \mathbb{K}\pi\otimes\mathbb{K}\pi \arrow[r, "\partial'_1"] \& \mathbb{K}\pi\otimes V\otimes\mathbb{K}\pi \arrow[r,"\partial'_0"] \& \mathbb{K}\pi \otimes\mathbb{K}\pi\arrow[r,"\mu"] \& \mathbb{K}\pi\arrow[r] \& 0\,,
	\end{tikzcd}
\end{align*}\vspace{-20pt}
\begin{align*}
	\partial'_1(1\otimes 1) &= \sum_{1\leq i\leq g} (a_1,b_1)\cdots(a_{i-1},b_{i-1})\left( 1\otimes a_i\otimes a_ib_ia_i^{-1}b_i^{-1} + a_i\otimes b_i \otimes b_ia_i^{-1}b_i^{-1}\right.\\
	&\hspace{100pt}\left.  - a_ib_ia_i^{-1}\otimes a_i \otimes b_i^{-1} - a_ib_ia_i^{-1}b_i^{-1}\otimes b_i \otimes 1\right)(a_{i+1},b_{i+1})\cdots(a_g,b_g),\\
	\partial'_0(1\otimes c\otimes 1) &= c\otimes c^{-1} -1\otimes 1 \quad \textrm{for } c\in\{a_i,b_i\}_{1\leq i\leq g}\,.
\end{align*}
Finally, swapping the second-to-left term by the automorphism
\begin{align*}
	\tau\colon \mathbb{K}\pi\otimes V\otimes\mathbb{K}\pi &\to \mathbb{K}\pi\otimes V\otimes\mathbb{K}\pi\\
	1\otimes c\otimes 1 &\mapsto 1\otimes c\otimes c^{-1} \quad \textrm{for } c\in\{a_i,b_i\}_{1\leq i\leq g}
\end{align*}
gives the proclaimed resolution: $d_1 = \tau\circ\partial'_1$ and $d_0 = -\partial'_0\circ\tau^{-1}$. \qed\\

Now we define a homological connection on $\Omega^1\mathbb{K}\pi$. Truncating (\ref{eq:resclean}) yields the resolution
\[
	\begin{tikzcd}[cramped,ampersand replacement = \&]
	0 \arrow[r] \& \mathbb{K}\pi\otimes\mathbb{K}\pi \arrow[r, "d_1"] \& \mathbb{K}\pi\otimes V\otimes \mathbb{K}\pi \arrow[r,"d_0"] \& \Omega^1\mathbb{K}\pi \arrow[r] \& 0\,.
	\end{tikzcd}
\]
On each degree, we set
\begin{align*}
	\nabla'_0(1\otimes a_i\otimes a_i^{-1}) &= 0,\; \nabla'_0(1\otimes b_i\otimes b_i^{-1}) = 0,\textrm{ and } \nabla'_1(1\otimes 1) = 0\,.
\end{align*}
$\nabla'_0$ is the push-out of $\nabla_\mathcal{C}$ by the natural map $p\colon\mathbb{K}F_{2g}\twoheadrightarrow\mathbb{K}\pi$. To compute the associated divergence, we need to lift Lie derivatives on $\Omega^1\mathbb{K}\pi$.

\begin{lemma}\label{lem:zetalift}
Let $f\in \Der_R(\mathbb{K}F_{2g})$, and $\hat f \in\Der_\mathbb{K}(\mathbb{K}\pi)$ the induced derivation. Then the $\hat f$-derivations defined by
\begin{align*}
	\lambda[\hat f]_0&\colon \mathbb{K}\pi\otimes V\otimes \mathbb{K}\pi\to \mathbb{K}\pi\otimes V\otimes \mathbb{K}\pi\colon 1\otimes v\otimes 1 \mapsto \sum_{c\in\{a_i,b_i\}_{1\leq i\leq g}} p(\partial'_cf(v))\otimes c \otimes p(\partial''_cf(v))\textrm{ and}\\
	\lambda[\hat f]_1&\colon \mathbb{K}\pi\otimes \mathbb{K}\pi\to \mathbb{K}\pi\otimes \mathbb{K}\pi\colon \hspace{58pt}1\otimes 1 \mapsto 0
\end{align*}
give a lift of $L_{\hat f}\colon \Omega^1\mathbb{K}\pi\to \Omega^1\mathbb{K}\pi$. 
\end{lemma}

\noindent Proof. First of all, $L_{\hat f}\circ d_0 = d_0\circ \lambda[\hat f]_0$ is clear from the construction. Next, we have
\begin{align*}
	(\lambda[\hat f]_0\circ d_1)(1\otimes 1) &= \lambda[\hat f]_0\left(\sum_{1\leq i\leq g} (a_1,b_1)\cdots(a_{i-1},b_{i-1}) \Big( 1\otimes a_i\otimes b_ia_i^{-1}b_i^{-1} + a_i\otimes b_i\otimes a_i^{-1}b_i^{-1}\right.\\
	&\left.\phantom{\sum_i}\hspace{50pt}- a_ib_ia_i^{-1}\otimes a_i\otimes a_i^{-1}b_i^{-1} - a_ib_ia_i^{-1}b_i^{-1}\otimes b_i\otimes b_i^{-1} \Big)(a_{i+1},b_{i+1})\cdots (a_g,b_g)\right)
\end{align*}
by the definition of $d_1$, and
\begin{align*}
	&= \sum_{1\leq i\leq g} \hat f((a_1,b_1)\cdots(a_{i-1},b_{i-1})) \Big( 1\otimes a_i\otimes b_ia_i^{-1}b_i^{-1} + a_i\otimes b_i\otimes a_i^{-1}b_i^{-1}\\
	&\phantom{\sum_i}\hspace{60pt}- a_ib_ia_i^{-1}\otimes a_i\otimes a_i^{-1}b_i^{-1} - a_ib_ia_i^{-1}b_i^{-1}\otimes b_i\otimes b_i^{-1} \Big)(a_{i+1},b_{i+1})\cdots (a_g,b_g)\\
	& \hspace{40pt} + (a_1,b_1)\cdots(a_{i-1},b_{i-1}) \lambda[\hat f]_0\Big( 1\otimes a_i\otimes b_ia_i^{-1}b_i^{-1} + a_i\otimes b_i\otimes a_i^{-1}b_i^{-1}\\
	&\phantom{\sum_i}\hspace{60pt}- a_ib_ia_i^{-1}\otimes a_i\otimes a_i^{-1}b_i^{-1} - a_ib_ia_i^{-1}b_i^{-1}\otimes b_i\otimes b_i^{-1} \Big)(a_{i+1},b_{i+1})\cdots (a_g,b_g)\\
	&\hspace{40pt} + (a_1,b_1)\cdots(a_{i-1},b_{i-1}) \Big( 1\otimes a_i\otimes b_ia_i^{-1}b_i^{-1} + a_i\otimes b_i\otimes a_i^{-1}b_i^{-1}\\
	&\phantom{\sum_i}\hspace{60pt}- a_ib_ia_i^{-1}\otimes a_i\otimes a_i^{-1}b_i^{-1} - a_ib_ia_i^{-1}b_i^{-1}\otimes b_i\otimes b_i^{-1} \Big)\hat f((a_{i+1},b_{i+1})\cdots (a_g,b_g))
\end{align*}
by the definition of $\lambda[\hat f]_0$, which is an $\hat f$-derivation. Using the isomorphism
\begin{align}\label{eq:dv}
	\Omega^1\mathbb{K}F_{2g} \cong \mathbb{K}F_{2g}\otimes V\otimes \mathbb{K}F_{2g}\colon dv\mapsto 1\otimes v\otimes 1\,,
\end{align}
the above equals
\begin{align*}
&\smash{\sum_{1\leq i\leq g}} \hat f((a_1,b_1)\cdots(a_{i-1},b_{i-1})) d(a_i,b_i) (a_{i+1},b_{i+1})\cdots (a_g,b_g)\\
&\quad\qquad + (a_1,b_1)\cdots(a_{i-1},b_{i-1}) df(a_i,b_i)(a_{i+1},b_{i+1})\cdots (a_g,b_g)\\
&\quad\qquad +  (a_1,b_1)\cdots(a_{i-1},b_{i-1})  d(a_i,b_i)\hat f((a_{i+1},b_{i+1})\cdots (a_g,b_g))\\
&= (p\otimes\id_V\otimes p)(df(\zeta)),
\end{align*}
which is zero since $f(\zeta)=0$. This concludes the proof.\qed\\

\begin{lemma}\label{lem:innerclosed}
Let $x\in\mathbb{K}F_{2g}$, $\hat x\in\mathbb{K}\pi$ its image and $\ad_{\hat x} = [\hat x\,,\cdot\,]$ the inner derivation by $\hat x$ on $\mathbb{K}\pi$. Then the $\ad_{\hat x}$-derivations defined by
\begin{align*}
	\lambda[\ad_{\hat x}]_0&\colon \mathbb{K}\pi\otimes V\otimes \mathbb{K}\pi\to \mathbb{K}\pi\otimes V\otimes \mathbb{K}\pi\colon 1\otimes v\otimes 1 \mapsto \sum_{c\in\{a_i,b_i\}_{1\leq i\leq g}} p(\partial'_c[x,v])\otimes c \otimes p(\partial''_c[x,v])\textrm{ and}\\
	\lambda[\ad_{\hat x}]_1&\colon \mathbb{K}\pi\otimes \mathbb{K}\pi\to \mathbb{K}\pi\otimes \mathbb{K}\pi\colon \hspace{58pt}1\otimes 1 \mapsto \hat x\otimes 1 - 1\otimes \hat x
\end{align*}
give a lift of $L_{\ad_{\hat x}}\colon \Omega^1\mathbb{K}\pi\to \Omega^1\mathbb{K}\pi$.
\end{lemma}

\noindent Proof. As in the previous lemma, $L_{\ad_{\hat x}}\circ d_0 = d_0\circ \lambda[\ad_{\hat x}]_0$ is clear. Next, we have
\begin{align*}
	(\lambda[\ad_{\hat x}]_0\circ d_1)(1\otimes 1) 
	&= \sum_{1\leq i\leq g} [\hat x,(a_1,b_1)\cdots(a_{i-1},b_{i-1})] \Big( 1\otimes a_i\otimes b_ia_i^{-1}b_i^{-1} + a_i\otimes b_i\otimes a_i^{-1}b_i^{-1}\\
	&\phantom{\sum_i}\hspace{40pt}- a_ib_ia_i^{-1}\otimes a_i\otimes a_i^{-1}b_i^{-1} - a_ib_ia_i^{-1}b_i^{-1}\otimes b_i\otimes b_i^{-1} \Big)(a_{i+1},b_{i+1})\cdots (a_g,b_g)\\
	& + (a_1,b_1)\cdots(a_{i-1},b_{i-1}) \lambda[\ad_{\hat x}]_0\Big( 1\otimes a_i\otimes b_ia_i^{-1}b_i^{-1} + a_i\otimes b_i\otimes a_i^{-1}b_i^{-1}\\
	&\phantom{\sum_i}\hspace{40pt}- a_ib_ia_i^{-1}\otimes a_i\otimes a_i^{-1}b_i^{-1} - a_ib_ia_i^{-1}b_i^{-1}\otimes b_i\otimes b_i^{-1} \Big)(a_{i+1},b_{i+1})\cdots (a_g,b_g)\\
	&+ (a_1,b_1)\cdots(a_{i-1},b_{i-1}) \Big( 1\otimes a_i\otimes b_ia_i^{-1}b_i^{-1} + a_i\otimes b_i\otimes a_i^{-1}b_i^{-1}\\
	&\phantom{\sum_i}\hspace{40pt}- a_ib_ia_i^{-1}\otimes a_i\otimes a_i^{-1}b_i^{-1} - a_ib_ia_i^{-1}b_i^{-1}\otimes b_i\otimes b_i^{-1} \Big)[\hat x,(a_{i+1},b_{i+1})\cdots (a_g,b_g)],\\
\intertext{which is equal to, using the isomorphism in (\ref{eq:dv}),}
	&\smash{\sum_{1\leq i\leq g}} [\hat x,(a_1,b_1)\cdots(a_{i-1},b_{i-1})] d(a_i,b_i) (a_{i+1},b_{i+1})\cdots (a_g,b_g)\\
	&\quad\qquad + (a_1,b_1)\cdots(a_{i-1},b_{i-1}) d[x, (a_i,b_i)](a_{i+1},b_{i+1})\cdots (a_g,b_g)\\
	&\quad\qquad +  (a_1,b_1)\cdots(a_{i-1},b_{i-1})  d(a_i,b_i)[\hat x,(a_{i+1},b_{i+1})\cdots (a_g,b_g)]\\
	=& \smash{\sum_{1\leq i\leq g}} \hat x(a_1,b_1)\cdots(a_{i-1},b_{i-1}) d(a_i,b_i) (a_{i+1},b_{i+1})\cdots (a_g,b_g)\\
	&\quad\qquad - (a_1,b_1)\cdots(a_{i-1},b_{i-1})\hat x d(a_i,b_i) (a_{i+1},b_{i+1})\cdots (a_g,b_g)\\
	&\quad\qquad + (a_1,b_1)\cdots(a_{i-1},b_{i-1}) dx(a_i,b_i)(a_{i+1},b_{i+1})\cdots (a_g,b_g)\\
	&\quad\qquad + (a_1,b_1)\cdots(a_{i-1},b_{i-1}) \hat x d (a_i,b_i)(a_{i+1},b_{i+1})\cdots (a_g,b_g)\\
	&\quad\qquad - (a_1,b_1)\cdots(a_{i-1},b_{i-1}) (a_i,b_i)dx (a_{i+1},b_{i+1}) \cdots (a_g,b_g)\\
	&\quad\qquad - (a_1,b_1)\cdots(a_{i-1},b_{i-1}) d (a_i,b_i)\hat x (a_{i+1},b_{i+1})\cdots (a_g,b_g)\\
	&\quad\qquad +  (a_1,b_1)\cdots(a_{i-1},b_{i-1})  d(a_i,b_i)\hat x (a_{i+1},b_{i+1})\cdots (a_g,b_g)\\
	&\quad\qquad -  (a_1,b_1)\cdots(a_{i-1},b_{i-1})  d(a_i,b_i)(a_{i+1},b_{i+1})\cdots (a_g,b_g) \hat x\\
	=& \smash{\sum_{1\leq i\leq g}} \hat x(a_1,b_1)\cdots(a_{i-1},b_{i-1}) d(a_i,b_i) (a_{i+1},b_{i+1})\cdots (a_g,b_g)\\
	&\quad\qquad -  (a_1,b_1)\cdots(a_{i-1},b_{i-1})  d(a_i,b_i)(a_{i+1},b_{i+1})\cdots (a_g,b_g) \hat x\\
	=& \;\hat xd_1(1\otimes 1) - d_1(1\otimes 1)\hat x\\
	=& \;(d_1\circ\lambda[\ad_{\hat x}]_1)(1\otimes 1).
\end{align*}
This concludes the proof.\qed\\

\begin{lemma}\label{lem:innerfree}
For a free-generating system $(x_i)_{1\leq i\leq 2g}$ of $F_{2g}$ and $y\in \mathbb{K}F_{2g}$, we have
\[
	\sum_{1\leq i\leq 2g}|\partial_i [y,x_i]| = (2g-1)|y\otimes 1 - 1\otimes y|
\]
in $|\mathbb{K}F_{2g}|^{\otimes 2}$.
\end{lemma}
\noindent Proof. Write $y = x_{i_1}^{\varepsilon_1}\cdots x_{i_r}^{\varepsilon_r}$ for some $\varepsilon_k\in\{1,-1\}$. Then, we have
\begin{align}\label{eq:innerdiv}
	\sum_{1\leq i\leq 2g}|\partial_i [y,x_i]| =  \sum_{1\leq i\leq 2g}|\partial_i(yx_i-x_iy)| = \sum_{1\leq i\leq 2g} |\partial_i(y)x_i - x_i \partial_i(y) + y\otimes 1 - 1\otimes y|.
\end{align}
For the first two terms on the right-hand side, we have
\begin{align*}
	\sum_{1\leq i\leq 2g}& |\partial_i(y)x_i - x_i \partial_i(y)|\\
	&= \sum_{1\leq i\leq 2g}| \partial_i (x_{i_1}^{\varepsilon_1}\cdots x_{i_r}^{\varepsilon_r})x_i - x_i \partial_i (x_{i_1}^{\varepsilon_1}\cdots x_{i_r}^{\varepsilon_r}) |\\
	&= \sum_{1\leq k\leq r} | x_{i_1}^{\varepsilon_1}\cdots \partial_{i_k}(x_{i_k}^{\varepsilon_k})\cdots x_{i_r}^{\varepsilon_r}x_{i_k} - x_{i_k} x_{i_1}^{\varepsilon_1}\cdots\partial_{i_k}(x_{i_k}^{\varepsilon_k})\cdots x_{i_r}^{\varepsilon_r} |
\end{align*}
since $\partial_i(x_{i_k}^{\varepsilon_k})$ is non-zero only if $i=i_k$. If $\varepsilon_k = 1$, the $k$-th term is equal to
\begin{align*}
	&|x_{i_1}^{\varepsilon_1}\cdots x_{i_{k-1}}^{\varepsilon_{k-1}}\otimes x_{i_{k+1}}^{\varepsilon_{k+1}}\cdots x_{i_r}^{\varepsilon_r}x_{i_k} - x_{i_k} x_{i_1}^{\varepsilon_1}\cdots x_{i_{k-1}}^{\varepsilon_{k-1}}\otimes x_{i_{k+1}}^{\varepsilon_{k+1}}\cdots x_{i_r}^{\varepsilon_r}|\\
	&\quad = |x_{i_1}^{\varepsilon_1}\cdots x_{i_{k-1}}^{\varepsilon_{k-1}}\otimes x_{i_k}x_{i_{k+1}}^{\varepsilon_{k+1}}\cdots x_{i_r}^{\varepsilon_r} - x_{i_1}^{\varepsilon_1}\cdots x_{i_{k-1}}^{\varepsilon_{k-1}} x_{i_k}\otimes x_{i_{k+1}}^{\varepsilon_{k+1}}\cdots x_{i_r}^{\varepsilon_r}|\,.
\end{align*}
If $\varepsilon_k = -1$, the $k$-th term is equal to
\begin{align*}
	&| - x_{i_1}^{\varepsilon_1}\cdots x_{i_{k-1}}^{\varepsilon_{k-1}}x_{i_k}^{-1}\otimes x_{i_k}^{-1}x_{i_{k+1}}^{\varepsilon_{k+1}}\cdots x_{i_r}^{\varepsilon_r}x_{i_k} + x_{i_k} x_{i_1}^{\varepsilon_1}\cdots x_{i_{k-1}}^{\varepsilon_{k-1}}x_{i_k}^{-1}\otimes x_{i_k}^{-1}x_{i_{k+1}}^{\varepsilon_{k+1}}\cdots x_{i_r}^{\varepsilon_r}|\\
	&\quad = | - x_{i_1}^{\varepsilon_1}\cdots x_{i_{k-1}}^{\varepsilon_{k-1}}x_{i_k}^{-1}\otimes x_{i_{k+1}}^{\varepsilon_{k+1}}\cdots x_{i_r}^{\varepsilon_r} + x_{i_1}^{\varepsilon_1}\cdots x_{i_{k-1}}^{\varepsilon_{k-1}} \otimes x_{i_k}^{-1}x_{i_{k+1}}^{\varepsilon_{k+1}}\cdots x_{i_r}^{\varepsilon_r}|\,.
\end{align*}
Therefore, in either way, they are equal to
\[
	|x_{i_1}^{\varepsilon_1}\cdots x_{i_{k-1}}^{\varepsilon_{k-1}}\otimes x_{i_{k}}^{\varepsilon_{k}}\cdots x_{i_r}^{\varepsilon_r} - x_{i_1}^{\varepsilon_1}\cdots x_{i_{k}}^{\varepsilon_{k}}\otimes x_{i_{k+1}}^{\varepsilon_{k+1}}\cdots x_{i_r}^{\varepsilon_r}|\,.
\]
by cancelling $x_{i_k}$ and $x_{i_k}^{-1}$ cyclically. Hence, the sum above is equal to 
\begin{align*}
	\sum_{1\leq k\leq r} &| x_{i_1}^{\varepsilon_1}\cdots x_{i_{k-1}}^{\varepsilon_{k-1}}\otimes x_{i_k}^{\varepsilon_k}\cdots x_{i_r}^{\varepsilon_r} - x_{i_1}^{\varepsilon_1}\cdots x_{i_k}^{\varepsilon_k}\otimes x_{i_{k+1}}^{\varepsilon_{k+1}}\cdots x_{i_r}^{\varepsilon_r} | = |1\otimes y - y\otimes 1|.
\end{align*}
Adding $2g|y\otimes 1 - 1\otimes y|$ in the third and fourth term in (\ref{eq:innerdiv}), we obtain the equality in the statement.\qed\\

\noindent\textbf{Proof of Theorem \ref{thm:cob}.} The left square is commutative by the construction of $\sigma$ and $v$. Next, by Lemmas \ref{lem:innerclosed} and \ref{lem:innerfree}, the divergence of an inner derivation takes its values in $|\mathbb{K}\pi|\otimes 1 + 1\otimes |\mathbb{K}\pi|$, so the map $\Div^{\nabla'}$ descends to $\mathrm{HH}^1(\mathbb{K}\pi)\to|\mathbb{K}\pi/\mathbb{K}1|^{\otimes 2}$. By Lemma \ref{lem:zetalift}, we have, for $f\in\Der_R(\mathbb{K}F_{2g})$,
\begin{align*}
	&\left(\lambda[\hat f]_1 - (i_{\varphi(\hat f)}\otimes \id)\circ\nabla'_1\right)(1\otimes 1) = 0\,,\\
	&\left(\lambda[\hat f]_0 - (i_{\varphi(\hat f)}\otimes \id)\circ\nabla'_0\right)(1\otimes c\otimes 1)\\
	&\qquad = \sum_{c'\in\{a_i,b_i\}_{1\leq i\leq g}} \Big(p(\partial'_{c'}f(c))\otimes c' \otimes p(\partial''_{c'}f(c))\Big) - 1\otimes c\otimes c^{-1}\hat f(c)\,,
\end{align*}
so that
\begin{align*}
	\Div^{\nabla'}(\hat f) = \sum_{c\in\{a_i,b_i\}_{1\leq i\leq g}} p(|\partial_cf(c)|) = p(\Div^{\nabla\!_\mathcal{C}}(f))
\end{align*}
modulo $|\mathbb{K}\pi|\otimes 1 + 1\otimes |\mathbb{K}\pi|$. Thus, the right square is also commutative; this completes the proof.\qed\\

\small
\bibliographystyle{alphaurl}
\bibliography{ncdtc.bib}

\end{document}